\documentclass[11pt]{amsart}
\usepackage{amsmath,amsthm,amscd,amssymb,amsfonts}
\usepackage{mathrsfs,upgreek}
\usepackage{hyperref}
\usepackage{mathtools}
\usepackage{graphicx}

\usepackage[sans]{dsfont} % Use \mathds{1} for indicator function. Erase/add [sans] for a different style
\usepackage[T1]{fontenc}
\DeclareMathAlphabet{\mathpzc}{OT1}{pzc}{m}{it} % Zapf Chancery math alphabet

\numberwithin{equation}{section} % Number formulas within sections
\numberwithin{figure}{section} % Number figures within sections

\theoremstyle{plain}
\newtheorem{theo}{Theorem}[section]
\newtheorem{prop}{Proposition}[section]

\newtheorem{lemm}[prop]{Lemma}

\newtheorem{theoalph}{Theorem}

\theoremstyle{definition}
\newtheorem{defi}[prop]{Definition}

\theoremstyle{remark}
\newtheorem{rema}[prop]{Remark}

\newtheoremstyle{citing}% name
  {3pt}%      Space above, empty = `usual value'
  {3pt}%      Space below
  {\itshape}% Body font
  {}%         Indent amount (empty = no indent, \parindent = para indent)
  {\bfseries}% Thm head font
  {.}%        Punctuation after thm head
  {.5em}%     Space after thm head: " " = normal interword space;
  {\thmnote{#3}}% Thm head spec

\theoremstyle{citing}
\newtheorem*{generic}{}% all text supplied in the note

\newcommand{\K}{\mathbb{K}}
\newcommand{\N}{\mathbb{N}}

\newcommand{\R}{\mathbb{R}}
\newcommand{\T}{\mathbb{T}}
\newcommand{\Z}{\mathbb{Z}}

\newcommand{\cB}{\mathcal{B}}

\newcommand{\cG}{\mathcal{G}}

\newcommand{\fG}{\mathfrak{G}}
\newcommand{\fH}{\mathfrak{H}}

\newcommand{\fL}{\mathfrak{L}}

\newcommand{\hC}{\widehat{C}}

\newcommand{\hE}{\widehat{E}}

\newcommand{\tL}{\widetilde{L}}

\newcommand{\teta}{\widetilde{\teta}}

\newcommand{\tLambda}{\widetilde{\Lambda}}

\newcommand{\tpi}{\widetilde{\pi}}

\renewcommand{\=}{\coloneqq}

\newcommand{\wtv}{\widetilde{v}}

\providecommand{\norm}[1]{\lVert#1\rVert}
\begin{document}
% \linenumbers

\title[Model sets with Euclidean internal space]{Model sets with Euclidean internal space}
\author{Mauricio Allendes Cerda}
\address{Mauricio Allendes Cerda, Departamento de Matem{\'a}ticas, Universidad Andres Bello, Avenida Rep\'ublica~498, Santiago, Chile}
\email{allendes.mauricio@gmail.com}
\author{Daniel Coronel}
\address{Facultad de Matem\'aticas, Pontificia Universidad Cat\'olica de Chile, Campus San Joaqu\'in, Avenida Vicu\~{n}a Mackenna 4860, 
Santiago, Chile.} 
\email{acoronel@mat.uc.cl}
%\address{Daniel Coronel, Departamento de Matem{\'a}ticas, Universidad Andres Bello, Avenida Rep\'ublica~498, Santiago, Chile}
%\email{daniel.coronel@mat-unab.cl}
%\author{}
%\address{}
%\email{}
% \thanks{}
% \urladdr{www}
% \date{\today}

\begin{abstract}
We  give an almost dynamical  characterization of inter-model sets with Euclidean internal space. This characterization is similar to previous results for general inter-model sets obtained independently by Baake, Lenz and Moody, and Aujogue. The new ingredients are two additional conditions on the rank of the Abelian group generated by the set of internal difference and  a flow on a torus defined via the   address map introduced by Lagarias that play the role of the maximal equicontinuous factor in the previous characterizations.
\end{abstract}
% \subjclass[2000]{S99}

\maketitle
% \setcounter{tocdepth}{1}
% \tableofcontents

%
%%% Main
%

\section{Introduction.}
In the  70's  Meyer introduced some separated nets in $\R^d$, also called Delone sets,    in connection with  his work in harmonic analysis. He  observed that each one of these sets, now called Meyer sets,  can be embedded into another type of Delone set called model set. This last collection is a sub-class of Meyer sets and are defined by  a simple geometric construction: they are the projection  on the first coordinate of some part of a lattice in $\R^d\times H$  where $H$, the internal space, is a locally compact Abelian group.

After the discovery of quasicrystals by D. Shechtman et al. \cite{Sh84}, model sets with  Euclidean internal space were proposed as geometric model for the atomic positions  in a quasicrystal.  Euclidean model sets and their associated dynamical systems played an important role in the mathematical diffraction theory of quasicrystals.   Hof in \cite{Hof} proved that every repetitive regular inter-model set (see definition in Sect. 2) has pure point diffraction and then, Schlottmann in \cite{Sch00} generalize this result  to repetitive regular   inter-model sets with arbitrary locally compact Abelian group as internal space.   
Euclidean model sets are also important in the theory of Pisot substitution  tilings. A central problem here has been  to understand when the  space generated by a Pisot substitution  is topologically conjugate to the space generated by a Euclidean model set \cite{BaKw06,BST10,ABBLS15}. 
%Recently, it has been proved that Euclidean model sets generated by windows with large boundaries can be very complex

In  \cite{Sch98}, Schlottmann gave a necessary and sufficient condition on a Delone set for being a general non-singular model set in terms of the recurrence structure of the  Delone set, and he asked for a characterization of non-singular model sets with well-behave internal space  as for example  $\R^n$.  
We recall that  every non-singular model set is a repetitive inter-model set.
 A dynamical characterization of repetitive regular inter-model sets was given by Baake, Moody and Lenz in \cite{BaLeMo} and then, Aujogue \cite{Au16a} extended this characterization to arbitrary repetitive inter-model sets not necessarily regular. Both results apply to general repetitive inter-model sets but left open the question of characterizing 
 repetitive inter-model sets with Euclidean internal space.
In this article, we answer this question by adding an algebraic and a dynamical property to the previous characterizations in \cite{BaLeMo} and \cite{Au16a}. The first condition is given in terms of the rank of the Abelian group generated by the set of difference of the Delone set and the second condition 
is  written in terms of a flow on a torus constructed from the 
address map introduced by Lagarias  in \cite{La99}, we call this flow the address system. We recall that every inter-model set is a Meyer set, and  all the previous characterizations of inter-model sets are written in the form of what we need to add to a Meyer set in order to have an inter-model set. Our result state that all the information needed for being an inter-model set with Euclidean internal space is encoded in the  rank of the group of differences and the dynamical relation  between the dynamical system associated to the Meyer set and the address system.

%Another motivation for knowing whether a Meyer set is an inter-model set with Euclidean internal space comes from the theory of substitution tilings and Delone sets. An important question here is to know when a substitution tiling or Delone set has pure point diffraction, in particular, in the case of Pisot substitutions. An approach to this problem is to prove that the   

In order to give a more detailed statement of our results we recall some definitions, see Sect. \ref{s:prelim} for details.

A discrete subset $\Lambda$ of $\R^d$ is a \emph{Delone set} if it is uniformly discrete and relatively dense. It is \emph{finitely generated} if the Abelian group generated by $\Lambda-\Lambda$ is finitely generated, and it is 
\emph{repetitive} if every pattern in $\Lambda$ appears with bounded gaps. Given a Delone set $\Lambda$ its \emph{hull} $\Omega_\Lambda$ is defined as the collection of all Delone sets whose local patterns agree with those of $\Lambda$ up to translation. If $\Lambda$ has \emph{finite local complexity}, then the hull can be endowed with a topology which is metrizable and compact. The subset of the hull of all Delone sets containing $0$ is called the \emph{canonical transversal} of $\Omega_\Lambda$ and we denote it by 
$\Xi_\Lambda$. The group $\R^d$ acts on the hull continuously by translation, given a (topological) dynamical system $(\Omega_\Lambda,\R^d)$. 
 Some combinatorial properties of the Delone set translate into dynamical properties. For example, repetitivity of $\Lambda$ is equivalent to minimality of $(\Omega_\Lambda,\R^d)$.
 It is well know in dynamical systems theory that, there is a dynamical system with an equicontinuous action of $\R^d$ that is a factor (semi-conjugacy) of $(\Omega_\Lambda,\R^d)$ and it is maximal with respect these properties.
 This dynamical system is unique up to topological conjugacy and we call it the \emph{maximal equicontinuous factor} of $(\Omega_\Lambda,\R^d)$. 

It is known that repetitivity implies finite local complexity, see for instance 
%\cite[Proposition $5.6$]{BaGr}
\cite{BaGr}, and that finite local complexity implies finitely generated, see \cite{La99}. 
A Delone set $\Lambda$ in $\R^d$ is a \emph{Meyer} set if
the set of difference $\Lambda-\Lambda$ is a Delone set.

Let $\Lambda$ be a finitely generated Delone set in $\R^d$. The \emph{rank} of $\Lambda$ is the rank of the Abelian group generated by $\Lambda$ as a subset of $\R^d$. We denote this group by $\langle \Lambda \rangle$, and by $s$ its rank. Let $\cB$ be a basis of $\langle \Lambda \rangle$, the \emph{address map}  for $\Lambda$ associated to $\cB$ is the coordinate map with respect to the basis $\cB$ from $\langle \Lambda \rangle$ to $\Z^s$.
 Notice that  since $\langle \Lambda \rangle$ is an Abelian group and $\langle \Lambda-\Lambda \rangle\subseteq\langle \Lambda \rangle$, we have that: $\langle \Lambda \rangle$ is finitely gerenated if and only if $\langle \Lambda-\Lambda \rangle$ is finitely generated. 
Also observe that if $\Lambda$ is a repetitive Meyer set in $\R^d$ then
 for every $\Lambda_0$ in $\Xi_\Lambda$ we have that 
$$\langle \Lambda_0 \rangle = \langle \Lambda_0 - \Lambda_0 \rangle = \langle \Lambda-\Lambda \rangle.$$  
Given a basis $\cB$ of $\langle \Lambda-\Lambda \rangle$, let $\varphi:\langle \Lambda-\Lambda \rangle\to \Z^s$ be the coordinate map with respect to the basis $\cB$.  We have that for every $\Lambda_0$ in $\Xi_\Lambda$ the address map of 
$\Lambda_0$ is equal to $\varphi$.

Lagarias proved in \cite{La99} that if $\Lambda$ is a Meyer set then, there is a linear map from $\R^d$ to $\R^s$ whose distance to the address map of  $\Lambda$ is uniformly bounded  on the points of $\Lambda$. Indeed, this property characterizes Meyer sets. Our first result gives the existence of one linear map that approximate the address map of 
all  Delone sets in $\Xi_\Lambda$, and it also gives a linear flow on a torus that we use to characterize inter-model sets with Euclidean internal space.

Put $\norm{x}_s$ for the Euclidean norm of $x$ in $\R^s$.

\begin{prop}[Address system]
\label{p:continuous eigenvalues}
Let $\Lambda$ be a repetitive Meyer set in $\R^d$ and let $s$ be the rank of  $\langle \Lambda-\Lambda \rangle$. 
Let $\cB$ be a basis of $\langle \Lambda-\Lambda \rangle$ and let $\varphi:\langle \Lambda-\Lambda \rangle\to \Z^s$ be the coordinate map with respect to the basis $\cB$.  
There are an injective linear map $\ell:\R^d\to \R^s$ and a constant $C>0$ such that for every $\Lambda_0$ in $\Xi_\Lambda$ and every $t\in\Lambda_0$ we have 
$$\|\varphi(t)-\ell(t)\|_s\le C.$$
Moreover, there is a linear flow $(\T^s,\R^d)$ defined by
$$(w,t)\in\T^s\times\R^d\longmapsto w+[\ell(t)]_{\Z^s},$$
and there is a homomorphism $\pi_{\text{Ad}}:\Omega_\Lambda \to \T^s$ such that for every $\Lambda'$ in $\Omega_\Lambda$ and every $t$ in $\R^d$ we have $\pi_{\text{Ad}}(\Lambda'-t) = \pi_{\text{Ad}}(\Lambda')+ [\ell(t)]_{\Z^s}.$
\end{prop}

Notice that the dynamical system $(\T^s,\R^d)$ and the homomorphism $\pi_{\text{Ad}}$ in Proposition \ref{p:continuous eigenvalues} depend on the basis $\cB$ chosen, however if we change the basis then the new system is topologically conjugate to the previous one.  We call any of these dynamical systems, the \emph{address system} of $\Lambda$, and to the map $\pi_{\text{Ad}}$ the \emph{address homomorphism} of $\Lambda$, which are well defined up to topological conjugacy. 
Observe that each coordinate of $\pi_{\text{Ad}}$ in Proposition \ref{p:continuous eigenvalues} gives a topological factor of $(\Omega_\Lambda, \R^d)$ onto the circle $\T$, however the address system is not necessarily a topological factor of $(\Omega_\Lambda, \R^d)$.  The minimality of $(\Omega_\Lambda, \R^d)$ implies that the address system of $\Lambda$ is a topological factor of $(\Omega_\Lambda,\R^d)$ if and only if the  address system of $\Lambda$ is minimal. Finally,
notice  that if we denote by $A$ the representative matrix of $\ell$ in the canonical basis and by $A^T$ the transpose then we have that the address system is minimal if and only if 
$\text{Ker}(A^T)\cap \Z^s = \{0\}$, which gives a simple way to check minimality of the address system.

The next theorem is the main result of the article, it characterizes inter-model sets with Euclidean internal space.

\begin{theoalph}
\label{t:model set}
 A repetitive Meyer set $\Lambda$ in $\R^d$ is an inter-model set with Euclidean internal space  if and only if   
 $\text{rank}(\langle\Lambda-\Lambda\rangle)>d$ and  the address system of $\Lambda$ is a topological factor of $(\Omega_\Lambda,\R^d)$ such that there is one point with a unique preimage under the factor map. In particular, if $\Lambda$ is an inter-model set with Euclidean internal space then the  address system of $\Lambda$ is the maximal equicontinuous factor of $(\Omega_\Lambda,\R^d)$.
 \end{theoalph}

%In \cite{Au16a}, Aujogue proved that a repetitive Meyer set is an inter-model set if and only if there is one point in its maximal equicontinuous factor with a unique preimages under the factor map. In Theorem \ref{t:model set}, we prove that the internal space is Euclidean with the additional conditions that $\text{rank}(\langle\Lambda-\Lambda\rangle)>d$ and the address system of $\Lambda$ is a topological factor of $(\Omega_\Lambda,\R^d)$. 
% In relation with the number of preimages of the address factor map  we have the following result.
% \begin{prop}
% \label{p:bounded preimages}
% Let $\Lambda$ be a repetitive Delone set  in $\R^d$. If  $\text{rank}(\langle\Lambda-\Lambda\rangle)>d$ and the address system of $\Lambda$  is a factor  of 
% $(\Omega_\Lambda,\R^d)$ then  there are a positive constant $C'$ and  a residual set $F$ in the address factor  of $\Lambda$ such that for every point in $F$ the number of preimages under the address factor  map is bounded by $C'$. 
%\end{prop} 

From Theorem \ref{t:model set} and \cite[Theorem $5$]{BaLeMo}, we obtain the following characterization 
for regular inter-model sets with Euclidean internal space. Observe that  if the address system of $\Lambda$ is minimal it is also uniquely ergodic.

\begin{theoalph}\label{t:regular model set}
  A repetitive Meyer set $\Lambda$ in $\R^d$ is a regular inter-model set with Euclidean internal space  if and only   
 $\text{rank}(\langle\Lambda-\Lambda\rangle)>d$ and  
  the address system of $\Lambda$ is a topological factor of $(\Omega_\Lambda,\R^d)$ such that the set of points in the address system with a unique 
  preimages under the factor map has full measure for the unique ergodic measure.
\end{theoalph}

For the proof of Theorem \ref{t:model set}, given a Meyer set we construct a cut and project scheme with a Euclidean internal space  and a window, which we call the ``Lagarias CPS'' and the ``minimal window'', respectively. What we actually prove in Theorem \ref{t:model set} is that if $\Lambda$ satisfies the necessary condition then it is an inter-model set generated by 
the Lagarias CPS and the minimal window. Using again \cite[Theorem $5$]{BaLeMo} we can give a more explicit version of Theorem \ref{t:regular model set}.

\begin{theoalph}\label{t:minimal window}
  A repetitive Meyer set $\Lambda$ in $\R^d$ is a regular inter-model set with Euclidean internal space  if and only   
 $\text{rank}(\langle\Lambda-\Lambda\rangle)>d$, the address system of $\Lambda$ is a topological factor of $(\Omega_\Lambda,\R^d)$ such that there is one point with a unique preimage under the factor map and the  boundary of minimal window of $\Lambda$ has measure zero.
\end{theoalph}
%This result can be seen as an extension of a similar result for Pisot substitutions. In this case it is known that the boundary of the Rauzy fractal has zero measure  
%\cite{}
%and the minimal  window is just the image of the Rauzy fractal  by a linear isomorphism \cite{}. Thus, for knowing whether
%

In order to put in context our results we mention an application to the theory of  unimodular Pisot substitution tilings.
In this setting one can prove that the address system corresponds to the canonical torus and it  is a topological factor of the hull of the substitution tiling. We can also prove that the minimal window is the image of the Rauzy fractal  by a linear isomorphism.
Using  Theorems \ref{t:model set}  and \ref{t:minimal window}  and the fact that the Rauzy fractal has zero measure boundary (see for instance \cite{BST10}) one can give another proof of  the following  known characterization of pure point unimodular Pisot substitution tilings as regular model sets with Euclidean internal space 
(see Theorem 7.3, Corollary 9.4, and Remark 18.6 in  \cite{BaKw06}).  

\begin{theo}\label{t:pisot}
 Let $\Omega_\Lambda$ be the hull of an  unimodular Pisot substitution tiling $\Lambda$ in $\R$. The following are equivalent:
\begin{itemize}
\item[(i)] $\Omega_\Lambda$ has pure point dynamical spectrum.

\item[(ii)] $\Omega_\Lambda$ is the hull of a regular model set with Euclidean internal space.

\item[(iii)] There is a point  in the address system  
of $\Lambda$ with a unique preimage under the factor map.
\end{itemize}
\end{theo}

\subsection{Strategy of the proof of Theorem \ref{t:model set} }
\label{ss:strategy}
The first step is to define the address system. This is done in Proposition \ref{p:continuous eigenvalues} where we showed that the address map on the canonical transversal minus the linear approximation define a continuous cocycle on the transverse groupoid. This cocycle is bounded. Then we apply a groupoid version of Gottschalk-Hedlund Theorem  \cite{Re12}  to prove that the cocycle is a coboundary. Using this coboundary we define the homomorphism.
For the necessary condition we show that  the maximal equicontinuous factor of the hull of an inter-model set with Euclidean internal space is topologically conjugate
to the address system.  
%The necessary condition in Theorem \ref{t:model set} is not difficult. 
To show that the condition is also sufficient is much more complicated. The first step is to prove that for a repetitive Meyer set the construction of an inter-model set given by Lagarias in \cite{La99}  gives a cut and project scheme (CPS) that we called the Lagarias CPS\footnote{Using the terminology of Lagarias we prove that repetitivity implies that the cut and project scheme is irreducible, see Sect.\ref{ss:lag cps}}, see definition in Sect. \ref{ss:CPS and MS}. Then, we show that the maximal equicontinuous factor of the hull associated to the Lagarias CPS (for any window) is topologically conjugate to the address system.
Then, we use \cite[Proposition 3.3]{Au16b} to prove that the closure of the projection  in the internal space of the lifting of the Meyer set to the product space is a window for 
the Lagarias CPS. Finally, elaborating on the ideas in \cite[Theorem 6.1]{Au16a} we show that if there is a point with a unique preimage under the maximal equicontinuous factor map of the hull of the Meyer set then the Meyer set is an inter-model set.

\subsection{Notes and references}
\label{ss:references}
From Proposition \ref{p:continuous eigenvalues} we have that for every repetitive Meyer set $\Lambda$ in $\R^d$ the dynamical system $(\Omega_\Lambda,\R^d)$ 
%has $s\ge d$ continuous eigenvalues in $\R^{d}$ corresponding to a basis of the internal difference group of $\Lambda$.
has $d$ continuous   linearly independent eigenvalues in $\R^{d}$. This   was first  proved by Kellendonk and Sadun in \cite{KeSa14} using pattern equivariant cohomological methods. Our proof relies on dynamical methods.

\subsection{Organization}
\label{ss:organization}
In Sect. \ref{s:prelim} we give some definitions and results about the theory of aperiodic order related to Delone sets and its associated dynamical systems. 
In Sect. \ref{s:proof theo eigenvalues} we prove Proposition \ref{p:continuous eigenvalues}. In Sect. \ref{ss:proof tms NC} we prove the necessary condition of 
Theorem \ref{t:model set}. 
In  Sect. \ref{ss:lag cps}, we describe the Lagarias CPS. The proof of the sufficient condition in Theorem \ref{t:model set} is in Sect. \ref{sss:Proof of sufficient condition} and it uses a result that we prove later in Sect. \ref{s:proof MTL}:  the Main Technical Lemma. 

%(Sect. \ref{l:factor}). In the last two sections we prove these results. In Sect. \ref{s:proof theo eigenvalues} we proved the Proposition \ref{p:continuous eigenvalues} using the Gottschalk-Hedlund's Theorem for groupoids \cite{Re12}. Finally, using the Lagarias CPS and some ideas from \cite{Au16b}, we proved the Main Technical Lemma in Sect. \ref{s:proof MTL}.

\subsection{Acknowledgments}
\label{ss:acknowledgments}
We would like to thank Marcy Barge for the references about substitution tilings.

The second named author acknowledges partial support from FONDECYT grant 1201612.

\section{Preliminaries.}
\label{s:prelim}
Let $\R^d$ be the Euclidean $d$-space  endowed with its Euclidean norm that we denote by $\|\cdot\|_d$.
\subsection{Delone sets}
\label{ss:delone}
A subset $\Lambda$ of $\mathbb{R}^d$  is called a  Delone set if it is \emph{uniformly discrete}, meaning that there is $r>0$
such that every closed ball of radius $r$ intersects $\Lambda$ in at most one point; and 
\emph{relatively dense}, which means that there is  $R>0$ such that 
every   closed ball of radius $R$ intersects $\Lambda$ in at least one
point. 

Let $\Lambda$ be a Delone set  in $\mathbb{R}^d$. For every $t\in \R^d$, 
we denote by $\Lambda-t$   the Delone set $\{x-t \mid x\in \Lambda\}$.
%We say that $\Lambda$ is \emph{aperiodic} 
%if it does not have translational symmetries, that is, if for every $t$ in $\R^d\setminus\{0\}$, we have $\Lambda \neq
%\Lambda-t$. 

For every $\rho>0$ and every $t$ in $\R^d$  denote by $B(t,\rho)$  the \emph{open ball} in $\mathbb{R}^d$ of radius $\rho$ and center $t$.  A \emph{$\rho$-patch} of $\Lambda$ centered at $t\in\R^d$ is the set $\Lambda\cap \overline{B(t,\rho)}$. 
We consider two notions of long-range order for Delone sets: 
The first one states that a Delone set $\Lambda$ has \emph{finite local complexity} if for every $\rho>0$ it has  a finite number of $\rho$-patches up to translation; 
and the second says that
$\Lambda$ is \emph{repetitive} if  for each $\rho>0$ there is a
number $M>0$ such that each
closed ball of radius $M$ contains the center of a translated copy of  every possible $\rho$-patch  of $\Lambda$. 
Observe that every repetitive Delone set has finite local complexity (see \cite[Proposition $5.6$]{BaGr}).

A Delone set $\Lambda$ is \emph{finitely generated} if the Abelian group generated by 
$\Lambda-\Lambda$ is finitely generated. We denote by $\langle \Lambda-\Lambda \rangle$ this group. Observe that
if $\langle \Lambda-\Lambda \rangle$  has finite rank the group $\langle \Lambda \rangle$ also has finite rank.
For every finitely generated Delone set we define its \emph{rank} as the rank of the group $\langle \Lambda \rangle$.
We recall the following proposition proved in \cite{La99}.
\begin{prop}
\label{p:flc}
Let $\Lambda$ be a Delone set.
\begin{enumerate}
\item $\Lambda$ has finite local complexity if and only if $\Lambda-\Lambda$ is discrete.
\item If $\Lambda$ has finite local complexity then  it is finitely generated.
\end{enumerate}
\end{prop}

\subsection{Meyer sets and address map}
\label{ss:meyer set}
Let $\Lambda$ be a Delone set in $\R^d$. We say that $\Lambda$ is a \emph{Meyer set} if there is a finite set $F$ in $\R^d$ such that  
$$
\Lambda-\Lambda\subseteq \Lambda + F.
$$ 

In \cite{Me72},  Meyer proved that every model set is a Meyer set.
Let  $\Lambda$ be a Meyer set in $\R^d$ with rank $s$, and let $\cB \=\{v_1,\ldots,v_s\}$ be a basis of $\langle \Lambda\rangle$. 
We recall that the address map for $\Lambda$ ( associated to the basis 
$\cB$) is the map $\varphi: \langle \Lambda\rangle \to \Z^s$ 
such that to every $x$ in $\langle \Lambda\rangle$ assigns its coordinates in the basis $\cB$. The following characterization of Meyer set  is used in the proofs of Proposition \ref{p:continuous eigenvalues} and Main Theorem.

\begin{theo}\cite[Theorem 3.1]{La99}
\label{teolag}
A Delone set $\Lambda$  in $\R^d$  is a Meyer set if and only if it finitely generated and every   address map 
$$
\varphi: \langle \Lambda \rangle \to \Z^s,
$$ 
is almost linear, that is,
there are a unique linear map $\ell:\R^d\to \R^s$ and a constant $C>0$ such that for every $x$ in $\Lambda$ we have
\begin{equation}\label{eq:condqualin}
\| \varphi(x) - \ell(x) \|_s \le C.
\end{equation}
\end{theo}
\begin{rema}
\label{r:ideal address}
In the proof of \cite[Theorem 3.1]{La99} it was proved that $\ell$ is some kind of ``ideal address map'' in the sense that if $\{v_1,\ldots, v_s\}$ is the basis of $\langle \Lambda \rangle$ that we used to define the address map of $\Lambda$  then  for every $t$ in $\R^d$ we have
 \begin{equation}
 \label{e:ideal address}
 \sum_{i=1}^s \ell_i(t) v_i = t.
 \end{equation}
\end{rema}

\subsection{Dynamical systems and transverse groupoid.}\label{ss:dysy-groupoid}
Let $\Lambda\subseteq\R^d$ be a Delone set with finite local complexity. The \emph{hull} of $\Lambda$ is the collection of all Delone sets in $\R^d$ whose $\rho$-patches, for every 
$\rho>0$, are also $\rho$-patches of  $\Lambda$ up to translation. We denote this set by $\Omega_\Lambda$. There is a natural metrizable topology on $\Omega_\Lambda$. Roughly speaking,  two Delone sets are close in this topology  if they agree on a large ball around the origen up to a small translation.
In particular, for every $\Lambda'$ in $\Omega_\Lambda$ a basis of  open neighborhoods for $\Lambda'$ is given by the following the sets:
First, for every $R>0$ put
$$
T(\Lambda',R)\=\{\tLambda \in \Omega_\Lambda \mid  \tLambda \cap \overline{B(0,R)} =   \Lambda' \cap \overline{B(0,R)} \},
$$
and for every $0<\varepsilon<R/2$ we define the open neighborhood  $N(\Lambda',\varepsilon, R)$ of $\Lambda'$ by
$$
N(\Lambda',\varepsilon, R) \=  \{\Lambda'' \in \Omega_\Lambda \mid \exists \tLambda \in T(\Lambda',R), \exists t\in B(0,\varepsilon), 
 \Lambda'' = \tLambda -t \},
$$  
for more details see for example  \cite{Sch00,FoHuKe,LeMo06,KeLe13}.  
If $\Lambda$ has finite local complexity then its hull $\Omega_\Lambda$ is compact. Observe  that  the action by translation of  $\R^d$ on $\Omega_\Lambda$ is continuous. Thus, we  obtain a topological dynamical system denote by $(\Omega_\Lambda,\R^d)$. The \emph{orbit} of   $x$ in $\Omega_\Lambda$ is the set $\{x-t \mid t\in \R^d\}$, and a subset $A$ of $\Omega_\Lambda$ is called \emph{invariant} if it is invariant by the action of $\R^d$. The dynamical systems $(\Omega_\Lambda,\R^d) $ is minimal if and only if the only closed invariant sets are the empty set and the whole space. It is well known that, minimality  is equivalent  to the fact every point has a dense orbit, and in the context of Delone sets repetitivity is equivalent to minimality. 

We recall that every topological dynamical system admits a \emph{maximal equicontinuous factor}. That is, a topological factor with an equicontinuous action such that any other equicontinuous factor is a topological factor of it (see for instance \cite{Ku03,BaKeSc,BaKe13}). For a topological dynamical system $(X,G)$ we denote by $(X_{\text{me}},G)$ its maximal equicontinuous factor. Given two minimal dynamical systems $(X,G)$ and $(Y,G)$, and a factor map $\pi:(X,G)\to (Y,G)$ we say that 
$\pi:(X,G)\to (Y,G)$ is an \emph{almost automorphic extension} or that   
 $(X,G)$ is an \emph{almost automorphic extension of} 
$(Y,G)$   if there is a point in  $Y$ with a unique preimage under $\pi$.

The \emph{transversal} of the hull is the closed subset 
$$\Xi_{\Lambda}:=\{x\in\Omega_\Lambda\mid 0\in x\}\subseteq\Omega_\Lambda.$$
In general, the restriction of the action of $\R^d$  to $\Xi_\Lambda$ is not defined. For this reason, to study the dynamical properties of the transversal we introduce the \emph{transverse groupoid},
$$\fG_\Lambda=\{(x,t)\in\Xi_{\Lambda}\times\R^d\mid x-t\in\Xi_{\Lambda}\}\subseteq\Xi_{\Lambda}\times\R^d.$$
This set,  endowed with the induced topology from the product space $\Xi_{\Lambda}\times\R^d$, has the  structure of a  topological groupoid 
(see \cite{Re80} for the abstract definition  of topological groupoids).  Two elements $(x,t)$ and $(z,s)$ in $\fG_\Lambda$ are composable if and only if $x-t=z$, and the composition of $(x,t)$ and $(z,s)$ is defined by
$$(x,t)\cdot(z,s)=(x,t+s).$$
The inverse map ${\cdot}^{-1}:\fG_\Lambda\to\fG_\Lambda$ is defined by $(x,t)^{-1}=(x-t,-t)$ 
 and the \emph{domain} $d:\fG_\Lambda\to\Xi_{\Lambda}$ and \emph{range} $r:\fG_\Lambda\to\Xi_{\Lambda}$ maps are defined by 
$$d(x,t)=x\quad\textrm{and}\quad r(x,t)=x-t.$$
Notice that  $d(\fG_\Lambda)=r(\fG_\Lambda)=\Xi_{\Lambda}$. In this context, the set $\Xi_{\Lambda}$ is called the \emph{unit space} of $\fG_\Lambda$. 

We say that a subset $E$ of the unit space  is \emph{invariant} by the groupoid $\fG$ if $E=r(d^{-1}(E))$. We recall the following definition from \cite{Re80}.
\begin{defi}
A groupoid is \emph{minimal} if the only open invariant subsets of its unit space are the empty set and the unit space itself.
\end{defi}

The following result relates the minimality of $(\Omega_\Lambda,\R^d)$ with the minimality of the transverse groupoid.

\begin{prop}\label{p:gr.min=dyn.min}
The topological groupoid $\fG_\Lambda$ is minimal if and only if the dynamical system $(\Omega_\Lambda,\R^d)$ is minimal.
\end{prop}
\proof 
First, observe that for every subset $E$ of $\Xi_\Lambda$ we have 
\begin{equation}\label{e:grupoid saturation}
r\left(d^{-1}(E)\right)=\{x-t\in\Xi_{\Lambda}\mid x\in E,t\in x\}.
\end{equation}
Assume that the dynamical system $(\Omega_\Lambda,\R^d)$ is minimal. Suppose, by contradiction, that $E\subseteq \Xi_\Lambda$  invariant by the grupoid 
$\fG_\Lambda$. 
Define 
$$
\hE\= \{ x-t  \in \Omega_\Lambda \mid x\in E, t\in \R\}.
$$
We have that $\hE$ is open in $\Omega_\Lambda$ and by \eqref{e:grupoid saturation} it is invariant for the $\R^d$-action on $\Omega_\Lambda$. Then, the complement of $\hE$ is an invariant  non empty  closed set  strictly contained in $\Omega_\Lambda$ which contradicts the minimality of   $(\Omega_\Lambda,\R^d)$. 

Reciprocally, suppose that $(\Omega_\Lambda,\R^d)$ is not minimal. Let $C\subseteq\Omega_\Lambda$ be an invariant non empty closed set strictly contained in 
$\Omega_\Lambda$. Put  $E=C^c\cap\Xi_\Lambda$. By \eqref{e:grupoid saturation} we have
$$E\subseteq r\left(d^{-1}(E)\right).$$
Since $C$ is invariant $\R^d$-action, we get $r\left(d^{-1}(E)\right)=E$. So, $E$ is a non empty open set strictly  contained in $\Xi_\Lambda$  invariant by the groupoid and thus $\fG_\Lambda$ is not minimal.
\endproof

\subsection{Cut and project scheme and inter-model sets}
\label{ss:CPS and MS}
A \emph{cut and project scheme} (CPS) over $\R^d$  is the data $(H,L)$ of a locally compact $\sigma$-compact Abelian group $H$, a discrete set $L\subseteq \R^d\times H$ 
with compact quotient $( \R^d\times H)/L$  
whose first coordinate projection on $\R^d$ is one-to-one and whose second coordinate projection on $H$ is dense. A compact subset $W$  of $H$ that is the closure of its interior is called a \emph{window} for the CPS. 
In the CPS the space $\R^d$ is called \emph{the physical space}, the locally compact Abelian group  $H$  is called
\emph{the internal space} and the set $L$ \emph{the lattice}. Following \cite{Au16b}, we have that  a CPS can also be described as a triple $(H,\Gamma,s_H)$ where $H$ is a locally compact $\sigma$-compact Abelian group, $\Gamma$ a countable subgroup of $\R^d$ and $s_{H}:\Gamma\to H$ a group morphism with range $s_{H}(\Gamma)$ dense in $H$ such that the graph
$$\cG(s_{H})\=\{(\gamma,s_{H}(\gamma))\in\R^d\times H\mid\gamma\in\Gamma\}$$
 is a \emph{lattice},  that is, a discrete and co-compact set.
 When $H$ is an Euclidean space $\R^n$, for some positive integer $n$, we say that $(H,\Gamma,s_{H})$ is an \emph{Euclidean CPS}.

Let $(H,\Gamma,s_{H})$ be a CPS with window $W$. For every $w$ in $H$, the projection on $\R^d$ of the set 
$\cG(s_H)\cap (\R^d\times (w+W))$ is called a \emph{model set}. 
More generally, for every subset $V$ of $H$ and every $w$ in $H$ denote by $\curlywedge(w+V)$ the set
$$
\curlywedge(w+V)\=\{t\in\Gamma\mid s_{H}(t)\in w+ V\}.
$$
\begin{defi}\label{d:intermodel set}
Let $(H,\Gamma,s_{H})$ be a CPS over $\R^d$  with window $W$. A Delone set $\Lambda\subseteq\R^d$ is called \emph{inter-model} set if exist $t\in\R^d$ and $w\in H$ such that
$$\curlywedge(w+\text{int}(W))-t\subseteq\Lambda\subseteq \curlywedge(w+W)-t.$$
\end{defi}

We say that an inter-model set $\Lambda$ is \emph{non-singular} or \emph{generic} if there is $(t,w)$ in $\R^d\times H$ such that 
$$\curlywedge(w+\text{int}(W))-t = \Lambda = \curlywedge(w+W)-t.$$
Observe that this is equivalent to the fact that the boundary of $w + W$ does not intersects the projection of $\cG(s_H)$ in $H$.
Additionally, if  the boundary of $w + W$ has zero Haar measure we say the inter-model set is \emph{regular}.
%Observe that all model set is an inter-model set, also its translations.
\begin{rema}\label{r:non-singular}
Notice that by Baire Theorem, the fact that $\partial W$ has empty interior and $s_H(\Gamma)$ is countable, the set 
$$
NS:=H \setminus \bigcup_{\gamma^*\in s_H(\Gamma)} \gamma^*-\partial W
$$
is a dense $G_\delta$-set   in $H$. Moreover, 
for every $w$ in $H$, the boundary of $w + W$ does not intersects the projection of $\cG(s_H)$ in $H$ if and only if $w\in NS$. In particular,
for every $(t,w)$ in $\R^d\times H$, the set~$\curlywedge(w+W)-t$ is a non-singular inter-model set if and only if $w\in NS$.
\end{rema}

The following two results are folklore.
\begin{prop}\label{p:unique hull G-MS}
Let $(\fH,\fL,s_{\fH})$ be a CPS over $\R^d$ with window $W$. The class of generic model sets generated by $(\fH,\fL,s_{\fH})$ and  window $W$ gives a unique hull, denoted by 
$\Omega_{MS}$. Every inter-model set generated by $(\fH,\fL,s_{\fH})$ and  the window $W$ is repetitive if and only if it belongs to $\Omega_{MS}$.  In particular, 
for every repetitive inter-model set $\Lambda$ generated by $(\fH,\fL,s_{\fH})$ and the window $W$ we have that $\Omega_\Lambda = \Omega_{MS}$, and
the dynamical system
$(\Omega_{MS},\R^d)$ is minimal.
\end{prop}

Let $(H,\Gamma,s_{H})$ be a CPS in $\R^d$  and consider the set $\T_{\cG}\=(\R^d\times H)/ \cG(s_{H})$ with an action of $\R^d$ given by translation on the first coordinate. 
Moreprecisely, for every  $s\in\R^d$ and every $[(t,w)]\in\T_{\cG}$ the action of $s$ on $[(t,w)]$ is
$$[(t,w)]\cdot s\=[(t,w)]+[(s,0)].$$
We say that $W'$ is \emph{irredundant} if the equation $W'+w=W'$ holds only for $w=0$ in $H$.

\begin{theo}\label{t:parametmodel}
Let $(\fH,\fL,s_{\fH})$ be a CPS over $\R^d$, let $W$ be an irrendundant window, and
let $\Omega_{MS}$ be the hull of the repetitive inter-model sets generated by  $(\fH,\fL,s_{\fH})$ and $W$. 
Then,  every point in $\Omega_{MS}$ is an inter-model set, and there exists a factor map $\pi:\Omega_{MS}\to\T_{\cG}$ such that  for every  $\Lambda'$ in $\Omega_{MS}$ there is $(t,w)$ in $\R^d\times H$ such that
 $\pi(\Lambda')=[(t,w)]$
 if and only if
\begin{equation}
\curlywedge(w+\text{int}(W)) - t \subseteq\Lambda'\subseteq \curlywedge(w+W) - t.
\end{equation}
Moreover, the map $\pi$ is injective precisely on the subset of non-singular inter-model sets in $\Omega_{MS}$ and the dynamical system  $(\T_{\cG},\R^d)$ is the maximal equicontinuous factor of $(\Omega_{MS},\R^d)$.
\end{theo}
The proof of the Theorem \ref{t:parametmodel} is mainly in \cite{Sch00}. The proof that $(\T_{\cG},\R^d)$ is the maximal equicontinuous factor of $(\Omega_{MS},\R^d)$  follows from  the fact 
 that $(\T_{\cG},\R^d)$ is an equicontinuous factor and from the existence of points where $\pi$ is injective (see for instance \cite[Lemma 3.11]{AuBaKeLe}).

\subsection{Torus parametrization}\label{ss:torus parametrization}
The notion of torus parametrization was introduced in \cite{BaHePl}.  
Here, we recall its definition and some properties  we  will use later.  
Let $X$ be a compact space and let $(X,\R^d)$ be a topological dynamical system under the action of $\R^d$ by the homeomorphisms $\{\rho_t\}_{t\in\R^d}$. Consider a compact Abelian group $\mathbb{K}$ with a minimal action of $\R^d$ by homeomorphisms $\{\kappa_t\}_{t\in\R^d}$. A \emph{torus parametrization} is a continuous map $\pi:X\to\mathbb{K}$ such that for all $t\in\R^d$ and $x\in X$ we have 
$$\kappa_{t}\circ\pi(x)=\pi\circ\rho_{t}(x).$$
For more details see \cite{BaLeMo,Sch00}. We recall the following lemma.
\begin{lemm}\cite[Lemma 1]{BaLeMo} 
If $\pi:X\to\K$ is a torus parametrization then $\pi$ is onto.
\end{lemm}
Let  $\pi:X\to\mathbb{K}$ be a torus parametrization. A \emph{section} of $\pi$ is a map $s:\mathbb{K} \to X$ such that $\pi\circ s$ is the identity on $\mathbb{K}$. A point $x\in X$ is called \emph{singular} if the fiber $\pi^{-1}(\pi(x))$ contains more than one element. Otherwise, $x\in X$ is \emph{called non-singular}. The set of non-singular points of $X$ for $\pi$ is denoted $R_\pi(X)$. The following proposition was proved in \cite{BaLeMo}.

\begin{prop}\cite[Proposition 3]{BaLeMo}
\label{p:cont.tor.par}
Let $\pi:X\to\mathbb{K}$ be a torus parametrization and let $s$ be a section of $\pi$. Then $s$ is continuous at all points of $\pi(R_\pi(X))$.
\end{prop}

\section{The Address system.}
\label{s:proof theo eigenvalues}
In this section we prove  Proposition \ref{p:continuous eigenvalues}. Given a repetitive Meyer set $\Lambda$ in $\R^d$,  we  start defining a continuous and bounded cocycle in the transverse  groupoid of $\Lambda$. We use a version of Gottschalk-Hedlund's theorem for gropupoids to show that this cocycle is a coboundary. 
We  use this cocycle and the map defining the coboundary to construct an equicontinuous dynamical system and homomorphism from $(\Omega_\Lambda,\R^d)$ into this equicontinuous system.
%This gives $s$ continuous eigenvalues   of the transverse grupoid, where $s$ is the rank of $\langle \Lambda-\Lambda \rangle$.
%We  uses these continuous eigenvalues and their associated eigenfunctions to construct an equicontinuous factor $(\T^s,\R^d)$ of $(\Omega_\Lambda,\R^d)$. 

\subsection{Defining a cocycle on the groupoid.}\label{ss:def cocycle}
Let $\Lambda\subseteq\R^d$ be a repetitive Meyer set. 
Let $\cB=\{v_1,\dots,v_s\}\subseteq\R^d$  be a basis for $\langle \Lambda-\Lambda \rangle$ and let $\varphi:\langle \Lambda-\Lambda \rangle \to \Z^d$ be the coordinate map 
with respect to the basis $\cB$. Recall that by the repetitivity of $\Lambda$ for every $x\in\Xi_\Lambda$ we have that 
$\langle x-x \rangle=\langle \Lambda-\Lambda \rangle$ and thus,
the address map of $x$  associated to $\cB$ is equal to $\varphi$. Note that for  all $t$ and  $t'$ in $\langle \Lambda - \Lambda \rangle$ we have
\begin{equation}\label{e:propadresmap}
\varphi(t+t')=\varphi(t)+\varphi(t').
\end{equation}
From Theorem \ref{teolag}, for every $x\in \Xi_{\Lambda}$ there is a unique  linear map $\ell_x:\R^d\to\R^s$ such that 
\begin{equation}
\label{e:linear approx}
\xi_x \= \sup_{t\in x} \|\varphi(t)-\ell_x(t)\|_{s} < +\infty.
\end{equation}
We define the maps $\Phi:\fG_\Lambda\to\Z^s$ and  $L:\fG_\Lambda \to \R^s$ as follows:   for every $(x,t)\in\fG_\Lambda$, 
$$\Phi(x,t):=\varphi(t) \text{ and } L(x,t)\= \ell_x(t).$$
The aim of this subsection is to show that $L-\Phi$ define a continuous cocycle on $\fG_\Lambda$. For this, we first prove that $L$ does not depend on the first coordinate. The proof of the continuity is at the end of the subsection.
\begin{prop}\label{L=cte}
% Let $\Lambda$ be a repetitive Meyer set in $\R^d$ such that $\text{rank}(\langle\Lambda-\Lambda\rangle)=s$. 
% %Let $\cB$ be a basis for $\langle\Lambda-\Lambda\rangle$, and for every $x$ in $\Xi_\Lambda$ let $\varphi_x$ be the address map associated to $\cB$ and let $\ell_x$ the unique linear map approximating 
% %$\varphi_x$. 
% Then,
% there is a linear map  $\ell:\R^d\to \R^s$ 
% such that for all $(x,t) \in \fG_\Lambda$ we have $L(x,t)=\ell(t).$
There is a linear map  $\ell:\R^d\to \R^s$ 
 such that for all $(x,t) \in \fG_\Lambda$ we have $L(x,t)=\ell(t).$
\end{prop}

The proof of this proposition is given at the end of this subsection after some  lemmas.

\begin{lemm}\label{l:delonedenso}
Let $\Lambda'$ be  a relatively dense  set in  $\R^d$. The set $\left\{\frac{t}{\|t\|_{d}}\mid t\in \Lambda'\right\}$ is dense in the boundary of the Euclidean unitary ball centered on the origin. 
In particular, for all linear map $T:\R^d\to\R^s$ we have that 
  $$\|T\|_{op}  = \sup_{t\in x} \left \|T\left(\frac{t}{\|t\|_{d}}\right)\right \|_{s},$$
where $\|\cdot\|_{op}$ is the operator norm.
\end{lemm}
  \proof 
  Put $D\=\left\{\frac{t}{\|t\|_{d}}\mid t\in \Lambda'\right\}$. By contradiction 
  suppose the set $D$ is not dense in the boundary of $B(0,1)$. So, there exists an open set in the relative topology which contains no elements of
$D$.
If we proyect this open set towards infinity, it generates a cone that contains Euclidean balls of size arbitrarily large and where there are no points of $\Lambda'$. This contradicts 
the fact that  $\Lambda'$ is relatively dense.
\endproof 

%First, we use this lemma to prove that the linear map $\ell_x$ that approximates $\varphi_x$ is unique.

%\begin{lemm}\label{l:linearmap unique}
%Let $\Lambda'$ be a Meyer set in $\R^d$ and let $\varphi$ be an address map for $\Lambda'$. There is a unique linear map that approximate $\varphi$ in the sense of \eqref{eq:condqualin}. 
%\end{lemm}
%\proof Let $s$ be the rank of $\Lambda'$.
%Assume that there are linear map  $\ell_1$ and $\ell_2$ from $\R^d$ to $\R^s$ and constants $C_1>0$ and $C_2>0$ such that $\ell_1$ and $C_1$, and $\ell_2$ and $C_2$ 
%approximate $\varphi$ in the sense of \eqref{eq:condqualin}.
%Then, for every $t\in \Lambda'$ we have 
%$$\|\ell_1(t)-\ell_2(t)\|_s\le\|\varphi_x(t)-\ell_1(t)\|_s+\|\varphi_x(t)-\ell_2(t)\|_s\le C_1+C_2.$$
%By Lemma \ref{l:delonedenso} we have  $\norm{\ell_1-\ell_2}_{op}=0$ and then,  $\ell_1=\ell_2$.
%\endproof

\begin{lemm}\label{tl=orb}
% Let $\Lambda$ be a repetitive Meyer set in $\R^d$ such that $\text{rank}(\langle\Lambda-\Lambda\rangle)=s$. 
% Let $\cB$ be a basis for $\langle\Lambda-\Lambda\rangle$ and let $\varphi: \langle \Lambda - \Lambda \rangle \to \Z^d$ be the coordinate map with respect to $\cB$. 
% For every $x$ in $\Xi_\Lambda$  let $\ell_x$ the unique linear map  given by Theorem \ref{teolag} that approximates 
% $\varphi$ on $x$.
% For all $(x,t)\in\fG_\Lambda$, we have 
% $\ell_x=\ell_{x-t}.$
 For all $(x,t)\in\fG_\Lambda$, we have 
 $\ell_x=\ell_{x-t}.$
  \end{lemm}
\proof Fix $(x,t)$ in $\fG_\Lambda$. Let $u\in\R^d$ be such that $u\in x-t$. In particular, $t+u\in x$.  By \eqref{e:linear approx} we have
  $$\|\varphi(u)-\ell_{x-t}(u)\|_{s}\le \xi_{x-t}\quad\textrm{and}\quad\|\varphi(u)-\ell_{x}(t+u)\|_{s}\le \xi_{x}.$$
Using these inequalities  and \eqref{e:propadresmap}, we get
\begin{multline*}
   \|\ell_x(t+u)-\ell_{x-t}(u)\|_{s}\le\|\varphi(t+u)-\ell_x(t+u)\|_{s}+\|\varphi(t+u)-\ell_{x-t}(u)\|_{s}\\
                              \le \xi_x+\|\varphi(t)+\varphi(u)-\ell_{x-t}(u)\|_{s}\\
                               \le \xi_x+\|\varphi(t)\|_{s}+\xi_{x-t}.
\end{multline*}
Dividing by $\norm{u}_{d}$ on both sides of this last inequality, we obtain
$$\left\|\ell_x\left(\frac{t}{\norm{u}_{d}}\right)+\ell_x\left(\frac{u}{\norm{u}_{d}}\right)-\ell_{x-t}\left(\frac{u}{\norm{u}_{d}}\right)\right\|_{s}\le \frac{\xi_x+\norm{\varphi(t)}_{s}+\xi_{x-t}}{\norm{u}_{d}}.$$
Taking the limit when $\norm{u}_{d}\to +\infty$ we have
 $$\lim_{\substack{\norm{u}_{d}\to +\infty\\ u\in x-t}}\left\|(\ell_x-\ell_{x-t})\left(\frac{u}{\norm{u}_{d}}\right)\right\|_{s}=0.$$
 This together with Lemma \ref{l:delonedenso} implies that $\norm{\ell_x-\ell_{x-t}}_{op}=0$, and thus, concludes the proof of the lemma. 
\endproof

%Finally, to prove the proposition, define the map $L:\fG_\Lambda\to\R^s$ for any $(x,t)\in\fG_\Lambda$ by 
%$$L(x,t)=\ell_x(t).$$
%We want to prove that this function is contant in the first coordinate. That is to say, for all $x,y\in\Xi_{\Lambda}$ the linear transformations  $\ell_x$ and $\ell_y$ are the same.

\proof[Proof of Proposition \ref{L=cte}] 
Fix $y$ in $\Xi_\Lambda$. We prove that for every $x$ in $\Xi_\Lambda$ we have $\ell_x=\ell_y$.
By \eqref{e:propadresmap}, \eqref{e:linear approx} and  Lemma \ref{tl=orb}, for $t'$ in $y$ we have 
\begin{eqnarray}
\label{e:xi bound xi}
\xi_{y-t'}& =&\sup_{t\in y-t'}\left\|\varphi(t)-\ell_{y-t'}(t)\right\|_s\nonumber\\
      &=& \sup_{t\in y-t'}\left\|\varphi(t)-\ell_{y}(t)\right\|_s\nonumber\\
      &=& \sup_{t+t'\in y}\left\|\varphi(t+t')-\varphi(t')-\ell_{y}(t)\right\|_s\\
      &=& \sup_{t+t'\in y}\left\|\varphi(t+t')-\varphi(t')-\ell_{y}(t+t'-t')\right\|_s\nonumber\\             
      &=& \sup_{t+t'\in y}\left\|\varphi(t+t')-\varphi(t')-\ell_{y}(t+t')+\ell_{y}(t')\right\|_s\nonumber\\
      &\le&  2\; \xi_y.\nonumber
\end{eqnarray}

Fix $x$ in $\Xi_\Lambda$. By minimality, 
there is  a sequence $(t_n)_{n\in \N}$ in $\R^d$ such that $y-t_n$ converges to $x$ in $\Xi_\Lambda$.
Fix $t\in x$ and consider $\epsilon>0$ such that  $\norm{t}\le\frac{1}{\epsilon}$. There is  $N\in\N$ such that for all $n>N$ we have
 $$(y-t_n)\cap \overline{B\left(0,\frac{1}{\epsilon}\right)}=x\cap \overline{B\left(0,\frac{1}{\epsilon}\right)}.$$
 In particular, for all $n>N$ we get $t\in y-t_n$. Then, using Lemma \ref{tl=orb} and \eqref{e:xi bound xi}, for every $t$ in $x$ we have
$$\left\|\varphi(t)-\ell_{y}(t)\right\|_s=\left\|\varphi(t)-\ell_{y-t_n}(t)\right\|_s\le 2\;\xi_y.$$
By uniqueness of the map $\ell_x$,  we conclude the proof of the proposition.
 
 \endproof

Let $H$ be  an Abelian group. A \emph{cocycle} on the topological groupoid $\fG_\Lambda$ with values in $H$ is a map $c:\fG_\Lambda\to H$ such that for all composable pairs $(x,t)$ and 
$(z,s)$ in $\fG_\Lambda$ one has
$$c((x,t)\cdot(z,s))=c((x,t))+c((z,s)).$$

%For composable pairs $(x,t),(z,s)\in\fG_\Lambda$ we have, by linearity of $L$, that
%$$L((x,t)\cdot(z,s))=L((x,t+s))=\ell(t+s)=\ell(t)+\ell(s)$$
%and by \eqref{e:propadresmap} we know that 
%$$\Phi((x,t)\cdot(z,s))=\Phi((x,t+s))=\varphi_x(t+s)=\varphi_x(t)+\varphi_{x-t}(s).$$
%Thus, $L-\Phi$ satisfies the equation
%\begin{eqnarray}
%(L-\Phi)((x,t)\cdot(z,s)) &=& \ell(t)-\varphi_x(t)+\ell(s)-\varphi_{x-t}(s)\nonumber\\
%                          &=& (L-\Phi)(x,t)+(L-\Phi)(z,s).\nonumber
%\end{eqnarray}
%We can conclude that the function $(L-\Phi):\fG_\Lambda\to\R^s$ defined for $(x,t)\in\fG_\Lambda$ by 
%$(L-\Phi)(x,t)=\ell(t)-\varphi_x(t)$, is a cocycle in $\fG_\Lambda$.

%\subsection{Continuity and boundedness of the cocycle.}\label{ss:cont and bound cocycle}

\begin{lemm}\label{l:cont of cocycle}
%Let $\Lambda$ be a repetitive Meyer set in $\R^d$. The map $L-\Phi$ is a continuous cocycle on $\fG_\Lambda$.
The map $L-\Phi$ is a continuous cocycle on $\fG_\Lambda$.
\end{lemm}
\proof
By \eqref{e:propadresmap} and Proposition \ref{L=cte} we have that $L-\Phi$ is a cocycle.
Now we  prove the continuity of $L-\Phi$. Consider a sequence $\{(x_n,t_n)\}_{n\in\N}$ in $\fG_\Lambda$ that converges to $(x,t)$
 in $\fG_\Lambda$. By definition of convergence in the groupoid, we have that $\{x_n\}_{n\in\N}\subseteq\Xi_\Lambda$ converges to $x\in\Xi_\Lambda$, and $\{t_n\}_{n\in\N}$ converges to $t$ in $\R^d$. 
 Let $\epsilon$ be a positive real number less than the uniformly discrete radius of $\Lambda$ such that 
 $\norm{t}_d<\frac{1}{2\epsilon}$. There is a positive integer $N$ such that for all $n\ge N$ we have
\begin{equation}\label{e:conv in G}
x_n\cap\overline{B\left(0,\frac{1}{\epsilon}\right)}=x\cap\overline{B\left(0,\frac{1}{\epsilon}\right)}, \norm{t_{n}-t}_d<\epsilon \text { and }
\norm{t_n}_d<\frac{1}{\epsilon}.
\end{equation}
By definition of the groupoid $\fG_\Lambda$,  for all $n$ in $\N$ we have that $t_n\in x_n$, and also $t\in x$. By \eqref{e:conv in G}, for every $n\ge N$ we get $t_n=t$.
Then, for every $n\ge N$ we have
$$L(t_n)=\ell(t)\quad\text{and}\quad\Phi(x_n,t_n)= \varphi(t_n)= \varphi(t)=\Phi(x,t),$$
which implies the continuity of $L-\Phi$.

\endproof

\subsection{Proof of Proposition \ref{p:continuous eigenvalues}.}
We use the following version of Gottschalk-Hedlund's Theorem, due to Jean Renault, to find continuous eigenvalues of $\fG_\Lambda$. This version is adapted to our context from   \cite[Theorem 1.4.10]{Re80} and it appears in \cite{Re12}.

\begin{theo}\label{got-hed:groupoid}
Let $G$ be a minimal topological groupoid with compact unit space $X$. For a continuous cocycle $c:G\to\R^d$ the following properties are equivalent:
\begin{enumerate}
\item There exists a continuous function $g:X\to\R^d$ such that 
$$c=g\circ r-g\circ d.$$
\item There exists $x\in X$ such that $c(d^{-1}(x))$ is relatively compact.
\item $c(G)$ is relatively compact.
\end{enumerate}
\end{theo}

\proof[Proof of Proposition \ref{p:continuous eigenvalues}]
Let $\Lambda\subseteq\R^d$ be a repetitive Meyer set. 
Let $\cB=\{v_1,\dots,v_s\}\subseteq\R^d$  be a basis for $\langle \Lambda-\Lambda \rangle$ and let $\varphi:\langle \Lambda-\Lambda \rangle \to \Z^d$ be the coordinate map 
with respect to the basis $\cB$. Let $L$ and $\Phi$ be as in Sect. \ref{ss:def cocycle}.
We check that $\fG_\Lambda$ and the cocycle $L-\Phi:\fG_\Lambda \to \R^s$  verify the hypotheses of Theorem  \ref{got-hed:groupoid}).
By Proposition \ref{p:gr.min=dyn.min} the groupoid is minimal.
By Lemma \ref{l:cont of cocycle}, the map $L-\Phi$ is a continuous cocycle.
Let  $\ell$ be the linear map given by Proposition \ref{L=cte}.
By \eqref{e:linear approx}, for every $x\in\Xi_{\Lambda}$  the set 
$$
(L-\Phi)(d^{-1}(x)) = \{\ell(t)-\varphi(t) \mid t\in x\} 
$$
is bounded.
By Theorem \ref{got-hed:groupoid}, there is a continuous map $F:\Xi_{\Lambda}\to\R^s$ such that 
for every $(x,t)$ in $\fG_\Lambda$ we have
\begin{equation}
\label{e:coboundary}
\ell(t)-\varphi(t) = L(x,t)-\Phi(x,t)=F\circ r(x,t)-F\circ d(x,t) = F(x-t)-F(x).
\end{equation}
Since $F$ is continuous and the space $\Xi_\Lambda$ is compact there is a constant $C>0$ such that the inequality in the first part of Proposition \ref{p:continuous eigenvalues} holds.

%Let $A$ be the representative matrix of $L$ in the canonical bases and  for every $j$ in $\{1,\ldots, s\}$ denote by $A_{(j,\cdot)}$ the $j$-row of $A$.
%Recall that  the image of $\Phi$  is in $\Z^s$.  Hence, writting \eqref{e:coboundary} by coordinates and taking exponential we get that for every $j$ in $\{1,\ldots, s\}$,
% $$
% \exp(2\pi i\; F_j(x-t))=\exp(2\pi i\;\langle A_{(j,\cdot)}, t\rangle)\;\exp(2\pi i \;F_j(x), 
%$$
%and thus, for every  $j\in\{1,\dots,s\}$ the row vector $A_{(j,\cdot)}\in\R^d$ is a continuous eigenvalue for the transverse groupoid $\fG_\Lambda$ with  continuous eigenfunction 
%$$
%f_j(x)=\exp(2\pi i\;F_j(x)).
%$$
%By Lemma \ref{transeig=cascoeig} all these continuous eigenvalues are also continuous eigenvalues of $(\Omega_{\Lambda},\R^d)$. 
%For every  $j\in\{1,\dots,s\}$, we denote also  by $f_j$  the extension of the eigenfunction to $\Omega_\Lambda$.

Now we check that  $\ell$ is injective. By contradiction suppose that the kernel of $\ell$ has dimension greater than one. Hence, there is an infinite  subset of $\Lambda$ such that  the address map is bounded on this infinite set, which gives  a contradiction. 
%Thus, the matrix $A$ has $d$ linearly independent rows, which implies that  $(\Omega_{\Lambda},\R^d)$ has at least $d$  linearly independent continuous eigenvalues.

Finally we construct the address system. Denote by $\T^s$ the torus $\R^s/\Z^s$. Since $\ell$ is linear the following map define an equicontinuous action of $\R^d$ on $\T^s$: 
$$(w,t)\in\T^s\times\R^d\longmapsto w+[\ell(t)]_{\Z^s}.$$
Now we define $\pi_{\text{Ad}}:\Omega_{\Lambda}\to\T^s$
as follows: For every $y\in\Omega_{\Lambda}$ there exist $x\in\Xi_{\Lambda}$ and $t\in\R^d$ such that $y=x-t$,  put
$$
\pi_{\text{Ad}}(y) \= [F(x)]_{\Z^s}+[\ell(t)]_{\Z^s}.
$$
We verify that  $\pi_{\text{Ad}}$ is well defined. Indeed, suppose that for $y\in\Omega_{\Lambda}$ there are 
$x_1,x_2\in\Xi_{\Lambda}$ and  $t_1,t_2\in\R^d$ such that $y=x_1-t_1=x_2-t_2$. Thus, $x_1=x_2-(t_2-t_1)$, and by \eqref{e:coboundary} we have that 
$$
F(x_1)=F(x_2)+\ell(t_2-t_1)-\varphi(t_2-t_1),
$$
which is equivalent to
$$
F(x_1)+\ell(t_1)=F(x_2)+\ell(t_2)-\varphi(t_2-t_1).
$$
Together with the fact that $\varphi(t_2-t_1)\in\Z^s$, this implies that  $\pi_{\text{Ad}}$ is well defined. 
Now we prove the continuity of $\pi_{\text{Ad}}$. Fix $y\in\Omega_\Lambda$ and suppose that $y=x-t$ for some $x\in\Xi_\Lambda$ and $t\in\R^d$. For every $y'$ close to $y$ there is $x'$ in $\Xi_\Lambda$ close to $x$ and there is $t'$ close to $t$ such that 
$y'=x'-t'$. By the continuity of $F$ and $\ell$, the map $\tpi_\text{Ad}$ defined in a sufficiently small neighborhood of $y$ by $\tpi_\text{Ad}(y')=F(x')+\ell(t')$ is continuous. By the continuity of the canonical projection of $\R^s$ onto $\T^s$ we conclude that $\pi_\text{Ad}$ is continuous at $y$.
It remains to check that for every $y$ in $\Omega_\Lambda$ and every $t$ in $\R^d$ we have $\pi_\text{Ad}(y-t) = \pi_\text{Ad}(y) + [\ell(t)]_{\Z^s}$. 
Fix $y$ in $\Omega_\Lambda$ and fix $t$ in $\R^d$. There are $x_1$ and $x_2$  in $\Xi_\Lambda$ and $t_1$ and $t_2$ in $\R^d$ such that $y=x_1-t_1$ and $y-t=x_2-t_2$.
Then, $x_2 = x_1-(t_1-t_2+t)$. Using this,  \eqref{e:coboundary} and the fact that $\varphi(t_1-t_2+t)\in\Z^s$ we get that
\begin{multline*}
\pi_\text{Ad}(y-t) = [F(x_2)]_{\Z^s}+[\ell(t_2)]_{\Z^s} \\
= [F(x_1) + \ell(t_1-t_2+t) -\varphi(t_1-t_2+t) ]_{\Z^s}+[\ell(t_2)]_{\Z^s} \\
= [F(x_1) + \ell(t_1)  ]_{\Z^s}+[\ell(t)]_{\Z^s} =  \pi_\text{Ad}(y) + [\ell(t)]_{\Z^s},
 \end{multline*}
 which concludes the proof of the proposition.
\endproof

\section{Proof of Theorem \ref{t:model set}.}
\label{s:proof theo model set}
In this section we prove Theorem \ref{t:model set}. First, we prove a characterization of the maximal equicontinuous factor for an Euclidean CPS and then, we prove the necessary condition. After this, we use the address map to construct an Euclidean CPS and that we use in the proof of  the sufficient condition. Finally, we prove the sufficient condition assuming the Main Technical Lemma. This lemma is stated in Sect.\ref{ss:proof tms SC}, and it is proved in Sect.\ref{s:proof MTL}. 

%We assume in this section that $\Lambda\subseteq\R^d$ is a repetitive Meyer set and denote by $(\Omega_\Lambda,\R^d)$ its associated  dynamical system. 

\subsection{Necessary condition}\label{ss:proof tms NC}
Let $\Lambda$ be an inter-model set for an Euclidean CPS  over $\R^d$ with internal space $\R^n$, lattice $L$ and window $W$. 
Denote by $\Omega_{MS}$ the hull of the generic set generated by these data.
Repetitivity of $\Lambda$ and Proposition \ref{p:unique hull G-MS} imply that $\Omega_{MS}=\Omega_\Lambda$. By \cite[Theorem 8.1]{Au16a}, the associated dynamical system $(\Omega_{MS},\R^d)$ is almost automorphic (see also \cite{Sch00,FoHuKe}). 
The remaining part of the proof of  the necessary condition follows directly  from  Proposition \ref{p:conj MEF} below.

\begin{prop}\label{p:conj MEF}
Let  $\Omega_{MS}$ be the hull of the repetitive inter-model sets generated  by an Euclidean cut and project scheme  $(\R^n, \Gamma,s_{\R^n})$ over $\R^d$ and  a window $W$.
Then, for every $\Lambda$ in $\Omega_{MS}$ we have that  the group $\langle \Lambda-\Lambda \rangle$  is equal to $\Gamma$ and  its  rank  is $d+n$. Moreover,  the maximal equicontinuous factor of $(\Omega_{MS},\R^d)$ is topologically conjugate to  the address system of $\Lambda$.
\end{prop}
\proof 
Denote by $p_1$ and by $p_2$ the ortogonal projections from $\R^d\times\R^n$ onto $\R^d$ and $\R^n$, respectively, and put $L\= \cG(s_{\R^n})$.
Fix $\Lambda$ in $\Omega_{MS}$.
By \cite[Proposition 2.6 (ii)]{Mo97} for every $w$ in $\R^n$ we have that 
$$
\langle \curlywedge(w+W) \rangle = L.
$$
In particular, $\langle \curlywedge(w+W) - \curlywedge(w+W) \rangle = L$. By Proposition \ref{p:unique hull G-MS} there is $w$ in $NS$ such that   
$\curlywedge(w+W)$ is in $\Omega_{MS}$ and thus, by repetitivity 
$$
 \langle\Lambda-\Lambda\rangle = \langle \curlywedge(w+W) - \curlywedge(w+W) \rangle = L.
$$

Now, we prove that the maximal equicontinuous factor of $(\Omega_{MS},\R^d)$ is topologically conjugate to  the address system of $\Lambda$.
Fix a basis $\cB=\{\wtv_1, \ldots, \wtv_s\}$ of  $L$. Let $\ell$ be the linear map given by the Proposition \ref{p:continuous eigenvalues} applied to $\Lambda$ with the basis 
$p_1(\cB)$ for $\Gamma$ and let $(\T^s,\R^d)$ be the address system.
Denote by $\psi:\R^s\to\R^d\times \R^n$  the linear isomorphism sending the canonical basis of $\R^s$ onto $\{\wtv_1, \ldots, \wtv_s\}$, i.e.
$$\psi(u_1,\dots,u_s)=u_{1}\wtv_1+\cdots+u_{s}\wtv_s.$$
By \eqref{e:ideal address} for every $t\in\R^d$ we have  
\begin{equation}\label{eq:proy_2}
p_1(\psi(\ell(t)))=t.
\end{equation}
Define the map $\Psi:\T^s\to\;\T_{\cG}$ by $\Psi([w]_{\Z^s})=[\psi(w)]_{L}$.
%, where $[w]=w+\Z^s$ and $[\psi(w)]=\psi(w)+L$. 
Note that $\Psi$ is an homeomorphism.  
By \eqref{eq:proy_2}, for all $t\in\R^d$ and $[w]\in\T^{s}$,  we have
\begin{multline*}
\Psi([w]_{\Z^s}+[\ell(t)]_{\Z^s} ) = \Psi([w+\ell(t)]_{\Z^s} ) \\ = [\psi(w + \ell(t))]_L   =  [\psi(w)]_L +[\psi(\ell(t))]_L 
\\ = [\psi(w)]_L +\left[\left(p_1(\psi(\ell(t))),p_2(\psi(\ell(t)))\right)\right]_L \\ = [\psi(w)]_L +\left[\left(t,p_2(\psi(\ell(t)))\right)\right]_L.
\end{multline*}
For proving that $\Psi$  conjugates the address system with  the maximal equicontinuous factor $(\R^d\times \R^n / L, \R^d)$, we need to show that for every $t\in\R^d$,
$$p_2(\psi(\ell(t)))=0.$$
By Remark \ref{r:non-singular}, Proposition \ref{p:unique hull G-MS} and the fact that the window $W$ has non-empty interior there is $w$ in $NS$ such that $0\in w+W$ and the set $\curlywedge(w+W)$ is in $\Omega_{MS}$.  
Put $\Lambda_0\= \curlywedge(w+W)$. We have that $\Lambda_0$ is in $\Xi_\Lambda$. Observe that  $\varphi$ is also the address map for $\Lambda_0$ associated to  the basis  $p_1(\cB)$.
By Proposition \ref{p:continuous eigenvalues} there is  a constant $\hC>0$ such that  for every $t\in\Lambda_0$ we have
\begin{equation*}
\label{e:bounded lift}
\|p_2(\psi(\varphi(t)))-p_2(\psi(\ell(t)))\|_d\le\hC.
\end{equation*}
Together with the fact that $p_2(\psi(\varphi(\Lambda_0)))=p_2(s_{\R^n}(\Lambda_0))\subseteq w+W$ this implies that
the map $p_2\circ \psi \circ \ell$ is uniformly bounded on $\Lambda_0$. Using that $\Lambda_0$ is relatively dense in $\R^d$ and  that 
$p_2\circ \psi \circ \ell$  is linear, we get that  $p_2(\psi(\ell(\R^d)))$ is bounded, which implies  that $p_2(\psi(\ell(\R^d)))=0$. We conclude that $(\T^s,\R^d)$ and $(\T_{\cG},\R^d)$ are topologically conjugated finishing the proof of the lemma.
\endproof

\subsection{Sufficient condition}\label{ss:proof tms SC}

\subsubsection{The Lagarias cut and project scheme}
\label{ss:lag cps}
Let $\Lambda$ be  a Meyer set in $\R^d$ and 
suppose that $\langle \Lambda-\Lambda \rangle$ has rank $s>d$.  
Let $\cB$ be a basis of $\langle \Lambda-\Lambda \rangle$ formed by vectors $\{v_1,\ldots,v_s\}\subseteq\R^d$ and 
let $\varphi:\langle \Lambda-\Lambda \rangle \to \Z^d$ be the coordinate map 
with respect to the basis $\cB$.
Fix $\Lambda_0$  in $\Xi_\Lambda$. Remember that since $0\in\Lambda_0$, we have 
$\langle \Lambda-\Lambda \rangle=\langle\Lambda_0-\Lambda_0\rangle=\langle\Lambda_0\rangle$ and that $\varphi$ is also the address map for $\Lambda_0$.   Let $\ell:\R^d\to \R^s$ be the linear map given by  Proposition \ref{p:continuous eigenvalues}.
 Define $\phi:\R^s\to\R^d$ by $\phi(u_1,\dots,u_s)=u_{1}v_1+\cdots+u_{s}v_s$. By \eqref{e:ideal address} for every $t\in\R^d$ we have 
 \begin{equation}\label{eq:lattproy1}
 \phi\circ\ell(t)=t.
 \end{equation}
In particular,
 \begin{equation}\label{eq:ker ell0}
 \text{Ker}(\ell)=\{0\}\qquad\textrm{and} \qquad\text{Im}(\phi)=\R^d.
 \end{equation}
Put~$n\=s-d$~and note that  the dimension of~$\text{Ker}(\phi)$ is~$n$. Let $\cB'  \= \{k_1,\ldots,k_n\}$ be an orthonormal basis for $\text{Ker}(\phi)$. Notice that for every $1\le j\le s$ we have that the vector $w_j\=\ell(v_j)-e_j$ belongs to $\text{Ker}(\phi)$, where $e_j$ is the $j$th-canonical coordinate vector.
For every $j\in\{1,\ldots, s\}$ denote by $(\alpha_{j,1},\ldots \alpha_{j,n})$ the coordinates of $w_j$ in the basis $\cB'$, and define for every $j\in\{1, \ldots, s\}$ the vectors 
$$
v_j^\star \= (\alpha_{j,1},\ldots \alpha_{j,n})^t \text{ and } \wtv_j \= (v_j, v_j^\star).
$$
In the proof of  \cite[Theorem 3.1]{La99}, Lagarias proved  that the set $\widetilde{\cB}\=\{\wtv_1\ldots, \wtv_s\}$ is $\Z$-linearly independent in $\R^d\times \R^n$ and its generates a full rank lattice. Denote by $\tL$ the lattice generated by 
$\widetilde{\cB}$. Denote by $p_1$ and  $p_2$ the orthogonal projections  of $ \R^d\times \R^n$ onto $\R^d$ and $\R^n$, respectively. By construction $p_1$ is injective on $\tL$ and its image is $\langle \Lambda-\Lambda \rangle$. 
Denote by $\psi:\R^s\to\R^d\times \R^n$  the linear isomorphism sending the canonical basis of $\R^s$ onto $\{\wtv_1, \ldots, \wtv_s\}$, i.e.
$$\psi(u_1,\dots,u_s)=u_{1}\wtv_1+\cdots+u_{s}\wtv_s.$$
In  the proof of  \cite[Theorem 3.1]{La99}, it was proved that for every $t$ in $\langle \Lambda-\Lambda \rangle$ we have
\begin{equation}
\label{e:internal projection}
\| p_2(\psi(\varphi(t))) \|_n = \| \varphi(t) - \ell(t) \|_s.
\end{equation}
 
\begin{lemm}\label{l:CPS lagarias}
Let $\Lambda$ be  a repetitive Meyer set in $\R^d$. If the address system of $\Lambda$ is a topological factor of $(\Omega_\Lambda,\R^d)$, then $p_2(\tL)$ is dense in $\R^n$.
\end{lemm}
\proof
The proof is by contradiction. Assume that $p_2(\tL)$ is not dense. Then, there is a non empty closed ball $V\subseteq\R^n$ such that $p_2(\tL)\cap V=\{\emptyset\}$. In particular, 
\begin{equation}
\label{e:non-dense}
\tL\cap (\R^d\times V)=\{\emptyset\}.
\end{equation}
By Proposition \ref{p:continuous eigenvalues}  and \eqref{e:internal projection} there is  a constant $\hC>0$ such that  for every $t\in\Lambda_0$ we have
\begin{equation*}
\label{e:bounded lift'}
\max\{\| p_2(\psi(\varphi(t))) \|_n, \|p_2(\psi(\varphi(t)))-p_2(\psi(\ell(t)))\|_n \} \le\hC.
\end{equation*}
Therefore the linear map $p_2\circ\psi\circ\ell$ is uniformly bounded on $\Lambda_0$, which is relatively dense. Then, 
for all $t\in\R^d$ we have
\begin{equation}
\label{e: zero projection}
p_2\circ\psi\circ\ell(t)=0.
\end{equation}
Consider the dynamical system defined on the space $(\R^d\times\R^n)/\tL$  with the  following $\R^d$-action:
for every $ t \in \R^d$ and  every  $w\in (\R^d\times\R^n)/\tL,$
$$ w \cdot t \=w+[(t,0)]_{\tL}.$$
Define the map  $\Psi:\T^s\to(\R^d\times\R^n)/\tL$   by $\Psi([w]_{\Z^s})=[\psi(w)]_{\tL}$ for every $[w]_{\Z^s}$ in $(\R^d\times\R^n)/\tL$.
By \eqref{e: zero projection}  the map $\Psi$ is a topological conjugacy between the address system of $\Lambda$ and  the dynamical system just defined $((\R^d\times\R^n)/\tL, \R^d)$.
Let $\pi_\text{Ad}$ be   the address homomorphism  defined in Proposition \ref{p:continuous eigenvalues}. 
Since we are assuming that $\pi_\text{Ad}$ is  a factor of $(\Omega_\Lambda,\R^d)$, we have that the map 
$\Psi\circ\pi_\text{Ad}$ is also a factor from $(\Omega_\Lambda,\R^d)$ to  $((\R^d\times\R^n)/\tL, \R^d)$. 
By the repetitivity of $\Lambda$ we have that $(\Omega_\Lambda,\R^d)$ is minimal and then, the factor $((\R^d\times\R^n)/\tL, \R^d)$ is also minimal.
But the set $[\R^d\times V]_{\tL}$ is  closed  and $\R^{d}$-invariant, and  by \eqref{e:non-dense},  it is strictly contained in $(\R^d\times\R^n)/\tL$, which is a contradiction 
to the minimality of $((\R^d\times\R^n)/\tL, \R^d)$.

\endproof

Put $s_{L} \= p_2\circ \psi \circ \varphi$ on $\langle \Lambda - \Lambda \rangle$.
By Lemma \ref{l:CPS lagarias}
if the address system of $\Lambda$ is a topological factor of $(\Omega_\Lambda,\R^d)$ the triple $(\R^n,\langle \Lambda-\Lambda \rangle,s_L)$ is a CPS  and
 we call it the \emph{Lagarias CPS for} $\Lambda$.  

Recall that a window is irredundant if its redundancies group is trivial (see Sect. \ref{ss:CPS and MS}). 
By compactness, every window in $\R^n$ is irredundant. 
By Theorem  \ref{teolag}  and \eqref{e:internal projection}, the set 
$\overline{s_L(\Lambda_0)}\subseteq\R^n$ is a compact. Together with
 Proposition \ref{p:window} we obtain  the following result.
\begin{lemm}\label{p:window cps lagarias}
Let $\Lambda$ be a repetitive Meyer set in $\R^d$ and let $\langle \Lambda-\Lambda \rangle$ be the subgroup of $\R^d$ generated by $\Lambda-\Lambda$. 
Put $n=\text{rank}(\langle \Lambda-\Lambda \rangle)-d$ and assume that $n>0$.  Also assume that the address system of $\Lambda$ is a topological factor of $(\Omega_\Lambda,\R^d)$.
Let  $(\R^n,\langle \Lambda-\Lambda \rangle,s_L)$ be the Lagarias CPS  for $\Lambda$. For every $\Lambda_0$ in $\Xi_\Lambda$ the set 
$\overline{s_L(\Lambda_0)}$ is an irredundant window.
\end{lemm}
From the proof of Lemma \ref{l:CPS lagarias} and by Lemma \ref{p:window cps lagarias} we obtain the following lemma.
\begin{lemm}\label{l:me factors}
Let $\Lambda$ be a repetitive Meyer set in $\R^d$ and let $\langle \Lambda-\Lambda \rangle$ be the subgroup of $\R^d$ generated by $\Lambda-\Lambda$. 
Put $n=\text{rank}(\langle \Lambda-\Lambda \rangle)-d$ and assume that $n>0$. 
 Also assume that the address system of $\Lambda$ is a topological factor of $(\Omega_\Lambda,\R^d)$.
Let  $(\R^n,\langle \Lambda-\Lambda \rangle,s_L)$ be the Lagarias CPS  for $\Lambda$. For every $\Lambda_0$ in $\Xi_\Lambda$, let $\Omega_{MS}$ be the hull of the generic inter-model sets
generated by  $(\R^n,\langle \Lambda-\Lambda \rangle,s_L)$ and the window $\overline{s_L(\Lambda_0)}$. Then the maximal equicontinuous factor of $(\Omega_{MS},\R^d)$ is topologically conjugated to
the address system of $\Lambda$.
\end{lemm}

\subsubsection{Proof of sufficient condition}
\label{sss:Proof of sufficient condition}

The main technical step in the proof of the sufficient condition is the following lemma that we state below. Its proof will be given  in Sect.\ref{s:proof MTL}.
\begin{generic}[Main Technical Lemma]
\label{l:factor}
Let $\Lambda\subseteq\R^d$ be a repetitive Meyer set and let $\Gamma$ be the subgroup of $\R^d$ generated by $\Lambda$. Let $(H',\Gamma, s_{H'})$ be a CPS and suppose that $W'=\overline{s_{H'}(\Lambda)}$ is a window. Let $\Omega_{MS}$ be the hull of the generic model sets generated by $(H',\Gamma,s_{H'})$ and $W'$. Then, there is a factor map 
$$
\widetilde{\pi}: \Omega_{\Lambda}\to \Omega_{MS,\text{me}},
$$
such that if  $(\Omega_{\Lambda},\R^d)$ is an almost automorphic extension of  $(\Omega_{MS,\text{me}},\R^d)$ for $\tpi$, then 
there are  $\Lambda_0$ in $\Omega_{\Lambda}$   and  a non-singular inter-model set $\Lambda_1$  in $\Omega_{MS}$ such that $\Lambda_0=\Lambda_1$.
\end{generic}

\proof[Proof of sufficient condition in Theorem \ref{t:model set}]
Let $\Lambda$ be a repetitive Meyer set in $\R^d$ and let $\langle \Lambda-\Lambda \rangle$ be the subgroup of $\R^d$ generated by $\Lambda-\Lambda$.  Assume that~$\text{rank}(\langle \Lambda-\Lambda \rangle)=s>d$,  that the address system is a topological factor 
$(\Omega_\Lambda,\R^d)$ and  that $(\Omega_{\Lambda},\R^d)$ is almost automorphic  extension of the address system. 

 Let $(\R^n,\langle \Lambda-\Lambda \rangle, s_{L})$ be the Lagarias CPS for $\Lambda$ where $n=s-d$.  Fix $\Lambda_*$ in $\Xi_\Lambda$ and recall that by 
 the repetitivity of $\Lambda$ we have that $\Omega_\Lambda=\Omega_{\Lambda_*}$.  By  Lemma \ref{p:window cps lagarias} 
the set $W'=\overline{s_L(\Lambda_*)}$ is an irredundant window.   Denote by $\Omega_{MS}$ the hull of generic inter-model sets generated by 
$(\R^n,\langle \Lambda-\Lambda \rangle,s_L)$ and $W'$.
By Lemma \ref{l:me factors} the maximal equicontinuous factor of  $(\Omega_{MS}, \R^d)$ is topologically conjugated to the address system of $\Lambda$ which  agrees with the address system of $\Lambda_*$ by Proposition \ref{p:continuous eigenvalues}.
By hypothesis, the dynamical system $(\Omega_{\Lambda_*},\R^d)$ is an almost automorphic extension of the address systems of $\Lambda_*$ and then, it is also an almost automorphic extension of   
$(\Omega_{MS,me}, \R^d)$. By the Main Technical Lemma applied to $\Lambda_*$ and $(\R^n,\langle \Lambda-\Lambda \rangle, s_L)$, 
there are  $\Lambda_0\in\Omega_{\Lambda_*}$ and $\Lambda_1\in\Omega_{MS}$ such that $\Lambda_0=\Lambda_1$. By  the minimality of 
$(\Omega_{\Lambda_*},\R^d)$ we have
that  $\Omega_{\Lambda_*}$ is equal to the hull of  $\Lambda_0$ which is equal to  the hull of the generic model sets generated by a Euclidean CPS. Since $W'$ is irredundat  and $\Omega_\Lambda=\Omega_{\Lambda_*}$ by Theorem \ref{t:parametmodel} we conclude that $\Lambda$ is an inter-model set generated by a CPS with Euclidean internal space, finishing the proof of the sufficient condition.

\endproof

\section{Proof of Main Technical Lemma.}
\label{s:proof MTL}
In this section we prove the Main Technical Lemma used in the proof of  Theorem \ref{t:model set}. 
Indeed, we prove a more detailed version of the Main Technical Lemma for future references.

\begin{generic}[Main Technical Lemma']
\label{l:factor general version}
Let $\Lambda\subseteq\R^d$ be a repetitive Meyer set and let $\Gamma$ the subgroup of $\R^d$ generated by $\Lambda$.
Let $(H',\Gamma,s_{H'})$ be a CPS and suppose that $W'=\overline{s_{H'}(\Lambda)}$ is a compact, irredundant window  in $H'$. 

Let $\Omega_{MS}$ be the hull of the generic inter-model sets for the CPS $(H',\Gamma,s_{H'})$ and window $W'$.
Let  $\pi_0$ be the maximal equicontinuous factor map from $\Omega_{MS}$ to $\Omega_{MS,\text{me}}$, and denote by $R_{\pi_0}(\Omega_{MS})$ the set of non-singular points in $\Omega_{MS}$ for $\pi_0$. 
Then, there is a factor map 
$$
\widetilde{\pi}: \Omega_{\Lambda}\to \Omega_{MS,\text{me}}.
$$
Put~$\Omega_{\Lambda}^0 \= \widetilde{\pi}^{-1}(\pi_0(R_{\pi_0}(\Omega_{MS})))$. There is a continuous map
$$
\pi_1: \Omega_{\Lambda}^0\to R_{\pi_0}(\Omega_{MS})
$$
such that for every $\Lambda_0\in\Omega_{\Lambda}^0$ we have 
$$\pi_1(\Lambda_0-t)=\pi_1(\Lambda_0)-t\quad\text{and}\quad \widetilde{\pi}(\Lambda_0)= \pi_0\circ \pi_1(\Lambda_0).$$ 
Moreover, for every $\Lambda_1$ in $ R_{\pi_0}(\Omega_{MS})$  we have
\begin{equation}\label{e:mtl}
\Lambda_1=  \bigcup_{\Lambda'\in \tpi^{-1}(\pi_0(\Lambda_1))} \Lambda'.
\end{equation}
In addition, if $\tpi: \Omega_{\Lambda}\to \Omega_{MS,\text{me}}$ is an  almost automorphic extension  then  
$$ 
\pi_0(R_{\pi_0}(\Omega_{MS}))\cap \tpi(R_{\tpi}(\Omega_\Lambda))
$$ 
is a residual set in $\Omega_{MS,\text{me}}$ and 
for every  $\Lambda_1$  in 
$ R_{\pi_0}(\Omega_{MS}) $ such that $\pi_0(\Lambda_1) \in \tpi(R_{\tpi}(\Omega_\Lambda))$ we have that $\Lambda_1$ is in $\Omega^0_{\Lambda}$.
\end{generic}

The proof of the lemma will be given in Sect. \ref{ss:proof} after recalling the definition of optimal CPS of a Meyer set introduced in \cite{Au16a}.

\subsection{The optimal CPS and the optimal window}
\label{ss:optimal CPS} 
Let $\Lambda$ be a repetitive Meyer set in $\R^d$ and let $\Gamma$ be the subgroup of $\R^d$ generated by $\Lambda$. Define $\Xi^{\Gamma}$ as the collection of all $\Lambda'\in\Omega_\Lambda$ having support into $\Gamma$,
$$
\Xi^{\Gamma} \= \{\Lambda' \in \Omega_\Lambda \mid  \Lambda'\subseteq \Gamma \}.
$$
Observe that $\Xi_{\Lambda}\subseteq\Xi^{\Gamma'}$.
We consider the \emph{combinatorial topology} on $\Omega_\Lambda$, which  is obtained from the distance
$$
\text{dist}(\Lambda',\Lambda'')=\left\{\frac{1}{R+1}\mid \Lambda' \cap \overline{B(0,R)} = \Lambda'' \cap \overline{B(0,R)}\right\}.
$$ 
The combinatorial topology is always strictly finer than the usual topology on $\Omega_{\Lambda}$ and on the transversal  
$\Xi_\Lambda$ both topologies coincide.
We endow $\Xi^{\Gamma'}$ with the combinatorial topology. 
We say that $\Lambda'$ and $\Lambda''$ in $\Omega_\Lambda$ are \emph{strongly regionally proximal}, denoted $\Lambda'\sim_\text{srp}\Lambda''$, if for each $R>0$ there are $\Lambda_1,\Lambda_2\in\Omega_\Lambda$ and $t\in\R^d$ such that
$$
\begin{array}{c}
\Lambda'\cap \overline{B(0,R)}=\Lambda_1\cap\overline{B(0,R)}\\
\Lambda''\cap\overline{B(0,R)}=\Lambda_2\cap\overline{B(0,R)}\\
(\Lambda_{1}-t)\cap\overline{B(0,R)}=(\Lambda_{2}-t)\cap\overline{B(0,R)}.
\end{array}
$$
Since $\Lambda$ is a repetitive Meyer set we have that the strongly proximal relation  is a closed $\R^d$-invariant equivalent relation on $\Omega_\Lambda$, and moreover, it agrees with the equicontinuous relation, see \cite{BaKe13}. In particular, the quotient $\Omega_\Lambda/\sim_{\text{srp}}$ gives the maximal equicontinuous factor.

In Proposition \ref{p:optimal CPS} below, we recall some results in \cite{Au16a} which allow to introduce the optimal CPS and optimal window for a Meyer set. More precisely, part (1) is deduced  by \cite[Proposition 4.4 and Lemma 4.5]{Au16a}, part (2) comes from  \cite[Proposition 6.1 and Definition 6.2]{Au16a} and finally, part (3) is in   \cite[Theorem 7.1]{Au16a}.
\begin{prop}\label{p:optimal CPS}
Let $\Lambda$ be a repetitive Meyer set in $\R^d$ and let $\Gamma$ the subgroup of $\R^d$ generated by $\Lambda$.
\begin{enumerate}
\item If $\Lambda'\in\Xi^{\Gamma}$ then its equivalence class $[\Lambda']_{srp}$ is contained into $\Xi^{\Gamma}$.
\item The set $H\=\Xi^{\Gamma}/\sim_{\text{srp}}$ with the quotient topology admits a structure of locally compact abelian group such that $[\Lambda]_{srp}$ is the neutro, the map $s_{H}:\Gamma \to H$ defined by $s_{H}(\gamma)=[\Lambda-\gamma]_\text{srp}$ is a group morphism and $\overline{s_H(\Gamma)} = H$. 
\end{enumerate}
\end{prop}
We remark that  Aujogue defined  $s_H$  in \cite{Au16a} with a sign plus instead of a minus as we did. So, some results that we use from  \cite{Au16a} and  \cite{Au16b} look slightly different since we need to do a correction in the sign. 
From Proposition \ref{p:optimal CPS},  the triple $(H,\Gamma,s_{H})$ is a CPS. Moreover, by \cite[Theorem 6.3]{Au16a}, the set 
$[\Xi_{\Lambda}]_{\text{srp}}$ is a window for $(H,\Gamma,s_{H})$. The CPS $(H,\Gamma,s_{H})$ and the window 
$[\Xi_{\Lambda}]_{\text{srp}}$ are called the \emph{optimal CPS} and the \emph{optimal window} for $\Lambda$, respectively.
Indeed, in \cite{Au16b}, the author proved that the model set that it defines,
$$\underline{\Lambda}\=\{\gamma\in\R^d\mid s_H(\gamma)\in [\Xi_{\Lambda}]_{\text{srp}}\},$$
satisfies that for every model set $M$ that includes $\Lambda$ we have $\Lambda\subseteq\underline{\Lambda}\subseteq M$.

Finally, we recall some results in \cite{Au16b} that we use in the proof of the Main Technical Lemma'. 
The first result allows to prove that a compact and irredundant set is a window.
\begin{prop}\cite[Proposition 3.3]{Au16b}
\label{p:window}
Let $\Lambda$ be a repetitive Meyer set in $\R^d$ and let $\Gamma$ the subgroup of $\R^d$ generated by $\Lambda$. Let $(H,\Gamma, s_{H})$ and $W$ be  optimal CPS and window for $\Lambda$, respectively. Suppose that $(H',\Gamma,s_{H'})$ is a CPS such that the closure $W'$ of the set $s_{H'}(\Lambda)$ is compact and irredundant in $H'$. Then, there is a continuous open and onto morphism
$$
\theta: H\to H' 
$$
such that $s_{H'} = \theta \circ s_H$ on $\Gamma.$ Moreover, the set $W'$ is a window in $H'$ and $W'=\theta([\Xi_{\Lambda}]_{\text{srp}})$.
\end{prop}

In the following result we  recall  the definition of a   map that we use   to construct the maps $\pi_1$ and $\widetilde{\pi}$ in the statement of  the
Main Technical Lemma'.
\begin{lemm}\cite[Lemma 3.4, 3.5, 3.6]{Au16b}
\label{l:omega}
Let $\Lambda$ be a repetitive Meyer set in $\R^d$ and let $\Gamma$ the subgroup of $\R^d$ generated by $\Lambda$. 
%Let $(H,\Gamma, s_{H})$ and $W$ be  optimal CPS and window for $\Lambda_0$, respectively. 
Suppose that $(H',\Gamma,s_{H'})$ is a CPS such that the closure $W'$ of the set $s_{H'}(\Lambda)$ is compact and irredundant in $H'$. 
We have that each $\Lambda'$ in $\Xi^\Gamma$ defines a unique element $w_\Lambda'$ through
$$
\{w_\Lambda'\} = \bigcap_{\gamma\in \Lambda'} s_{H'}(\gamma)-W'.
$$
Define the map 
$$\begin{array}{cl}
\omega: & \Xi^\Gamma\to H'\\
        & \Lambda'\mapsto w_{\Lambda'}.
\end{array}$$
We have that $\omega$ is is  uniformly continuous for the combinatorial topology, and for all $\Lambda'\in\Xi^\Gamma$ and $\gamma\in\Gamma$ we have
\begin{enumerate}
\item $\omega(\Lambda'-\gamma)=\omega(\Lambda')-s_{H'}(\gamma).$
\item $\omega(\Lambda')=-\theta([\Lambda']_{\text{srp}})$, where $\theta$ is the morphism in Proposition \ref{p:window}.
\end{enumerate}
\end{lemm}

\subsection{Proof of Main Technical Lemma'}
\label{ss:proof}

Let $\Lambda\subseteq\R^d$ be a repetitive Meyer set and let $\Gamma$ be the subgroup of $\R^d$ generated by $\Lambda$. Let $(H',\Gamma,s_{H'})$ be a CPS and assume that $W'=\overline{s_{H'}(\Lambda)}$ is  a compact and irredundant window  in $H'$. 
Let $\Omega_{MS}$ be the hull of inter-model sets generated by $(H',\Gamma,s_{H'})$ and  $W'$. Recall that  the maximal equicontinuous factor $\Omega_{MS,\text{me}}$ 
can be obtained by the quotient 
 $(\R^d\times H')/ \cG(s_{H'})$ and denote 
by $\pi_0$ be the maximal equicontinuous factor map from $\Omega_{MS}$ to $\Omega_{MS,\text{me}}$.

\subsubsection{Construction of $\widetilde{\pi}$}
Now we construct the map $\widetilde{\pi}: \Omega_\Lambda \to \Omega_{MS,\text{me}}$.
For every $(t,w)$ in $\R^d\times H'$ we denote by $[(t,w)]$ its  equivalent class in   $\Omega_{MS,\text{me}}$.
For every $\tLambda$ in $\Omega_{\Lambda}$ there is $t\in\R^d$ such that $\tLambda-t$ is in $\Xi^\Gamma$,  define $\widetilde{\pi}(\tLambda)$ by
$$\widetilde{\pi}(\tLambda) \= [(-t,\omega(\tLambda-t))]\in\Omega_{MS,\text{me}}.$$
We verify  that $\widetilde{\pi}$ is well defined. Assume that there is $s$ in $\R^d$ such that $\tLambda-s$ is in $\Xi^\Gamma$. 
Observe that $t-s$ is in $\Gamma$. By part (1) in Lemma \ref{l:omega}, we have that 
\begin{multline*}
(-t,\omega(\tLambda-t)) = (-t+s-s,\omega(\tLambda-(t+s-s)))\\
 = (-s - (t-s),\omega(\tLambda-s)-s_{H'}(t-s)) \\
= (-s,\omega(\tLambda-s)) - ( t-s,s_{H'}(t-s)).
\end{multline*}
Since $( t-s,s_{H'}(t-s))$ belongs to $\cG(s_{H'})$, we have that 
$$[(-t,\omega(\tLambda-t))] = [(-s,\omega(\tLambda-s))],$$
and hence $\widetilde{\pi}$ is well defined. 

Now we check that $\widetilde{\pi}$ commutes with the $\R^d$ action on $\Omega_{\Lambda}$ and on $\Omega_{MS,\text{me}}.$
Let $\tLambda$ be in  $\Omega_{\Lambda}$ and  $t$ be in $\R$. There are $s$ and $s'$ in $\R^d$ such that $\tLambda-s$ 
and  $(\tLambda-t)-s'= \tLambda-(t+s')$ are in $\Xi^\Gamma$. 
Notice that $t+s'-s$ belongs to $\Gamma$.
Again, by part (1) in Lemma \ref{l:omega} we have
\begin{multline*}
(-s',\omega((\tLambda-t)-s')) = (-s',\omega((\tLambda-s)-(t+s'-s)))\\
= (-s',\omega(\tLambda-s)-s_{H'}(t+s'-s))\\
= (-s'+(t+s'-s),\omega(\tLambda-s))-(t+s'-s, s_{H'}(t+s'-s))\\
= (t-s,\omega(\tLambda-s))-(t+s'-s, s_{H'}(t+s'-s)).
\end{multline*}
Since $(t+s'-s, s_{H'}(t+s'-s))$ is in $\cG(s_{H'})$ we have 
$$ \widetilde{\pi}(\tLambda-t) = [(-s',\omega((\tLambda-t)-s'))] = [(-s,\omega(\tLambda-s))]+[(t,0)] = \widetilde{\pi}(\tLambda) + [(t,0)].$$
Now we prove that $\widetilde{\pi}$ is continuous. Let  $\Lambda'$ be $\Omega_{MS}$  and let $U$ be a neighborhood of $0$ in $\Omega_{MS,\text{me}}$. We can assume that $U=[B(0,r_0)\times U_{H'}]$ where $r_0>0$ and $U_{H'}$ is a neighborhood of $0$ in $H'$. 
There exists $t'\in\Lambda'$ such that $\Lambda'-t'\in\Xi_{\Lambda}\subseteq \Xi^\Gamma$. For $r>0$, denote 
$$C_r=\bigcap_{\gamma\in (\Lambda'-t')\cap\overline{B(0,r)}}s_{H'}(\gamma)-W',$$
and observe that for $r>r'$ we have $C_r\subseteq C_{r'}$. By Lemma \ref{l:omega},
\begin{equation}\label{e:transversal intersection}
\bigcap_{r>0}C_r=\{\omega(\Lambda'-t')\}.
\end{equation}
Now we prove that there is $r'>0$ such that for every $r\ge r'$ 
\begin{equation}\label{e:transversal continuity}
C_r \subseteq \omega(\Lambda'-t')+U_{H'}.
\end{equation}
By contradiction suppose that there is an increasing sequence $(r_i)_{i\in\N}$ of positive real numbers converging to infinity as $i$ goes to infinity, such that $(C_{r_i}-\omega(\Lambda'-t'))\cap U_{H'}^c\neq\emptyset$. Then, for every $i\in\N$ there is 
$$
x_{i}\in (C_{r_i}-\omega(\Lambda'-t'))\cap U_{H'}^c.
$$ 
Since for every  $i,j$ in $\N$  with $j\ge i$ we have $C_{r_j}\subseteq C_{r_i}$. 
By compactness of $C_{r_1}$ there is an accumulation point  $\tilde{x}$ of $(x_{i})_{i\in\N}$ in $U_{H'}^c$ and thus, $\tilde{x}\neq 0$. But $\tilde{x}$ also belongs to 
$\bigcap_{r>0} C_r -\omega(\Lambda'-t')$ which is $\{0\}$ by \eqref{e:transversal intersection}, giving the desired contradiction.

Put $R\= \norm{t'}_d+r' + r_0$ and consider 
set 
$$
T\= \{\tLambda \in \Omega_\Lambda \mid \Lambda' \cap \overline{B(0,R)} =  \tLambda \cap \overline{B(0,R)}\}.
$$
For every $\varepsilon>0$ sufficiently small the set 
$$
V_\varepsilon \=\{\Lambda'' \in \Omega_\Lambda \mid  \exists \tLambda \in T, \exists t \in B(0,\varepsilon), \Lambda'' = \tLambda-t\}
$$
is an open neighborhood of $\Lambda'$. Fix $\varepsilon < r_0$. By the definition of $R$, for every $\Lambda'' $ in $V_\varepsilon$ there are $t$ in $B(0,\varepsilon)$ and $\tLambda$ in $T$ such that  
  $$
  (\Lambda''-(t'-t)) \cap \overline{B(0,r')} = (\tLambda-t' )\cap \overline{B(0,r')} =  (\Lambda'-t' )\cap \overline{B(0,r')}.
  $$
Put  $t''\= t'-t$, we have $\norm{t'-t''}_d<r_0$ and  since $\Lambda'-t'$ is in $\Xi_\Lambda$ we also have that $\Lambda''-t''$ is in $\Xi_\Lambda \subseteq \Xi^{\Gamma}$.
Then,
$$\bigcap_{\gamma\in (\Lambda''-t'')\cap\overline{B(0,r')}}s_{H'}(\gamma)-W'=\bigcap_{\gamma\in (\Lambda'-t')\cap\overline{B(0,r')}}s_{H'}(\gamma)-W'.$$
Together with \eqref{e:transversal continuity}  this implies $\omega(\Lambda''-t'')\in\omega(\Lambda'-t')+U_{H'}$. Therefore,
 $\widetilde{\pi}(\Lambda'')=[-t'',\omega(\Lambda''-t'')]$ is included in
\begin{multline*}
[-t'+(t'-t''),\omega(\Lambda'-t')+U_{H'}]\subseteq[-t'+B(0,\delta),\omega(\Lambda'-t')+U_{H'}]\\
=[-t',\omega(\Lambda'-t')]+[B(0,\delta),U_{H'}],
\end{multline*}
showing the continuity of $\widetilde{\pi}$ at $\Lambda'$ in $\Omega_{MS}$.
 
 Finally, since the $\R^d$-action on $\Omega_{MS,\text{me}}$  is minimal  we have that $\widetilde{\pi}$ is surjective, which concludes the proof that $\widetilde{\pi}$ is a factor map.

\subsubsection{Definition of $\pi_1$}
Recall that $R(\Omega_{MS})$ denotes the set of non-singular points of $\Omega_{MS}$ for $\pi_0$ as defined in Sect.\ref{ss:torus parametrization}. By definition, all sections of $\pi_0$ agree on $\pi_0(R(\Omega_{MS}))$. Let $\widetilde{s}:\Omega_{MS,\text{me}} \to \Omega_{MS}$ be a section of $\pi_0$. Put $\Omega_{\Lambda}^0 \= \widetilde{\pi}^{-1}(\pi_0(R(\Omega_{MS}))$, and define the surjective map $\pi_1: \Omega_{\Lambda}^0\to R(\Omega_{MS})$ by $\pi_1 \= \widetilde{s} \circ \widetilde{\pi}$.

By the continuity of  $\widetilde{\pi}$ and  Proposition \ref{p:cont.tor.par} the map $\pi_1$ is also continuous.
Since $\widetilde{s}$ is a section of $\pi_0$, for every $\Lambda'$ in $\Omega_{\Lambda}^0$ we have
\begin{equation}\label{eq: rel pi_0 pi_1}
\widetilde{\pi}(\Lambda') = \pi_0 \circ \pi_1(\Lambda'). 
\end{equation}
Since $\widetilde{s}$ commutes with the action of $\R^d$ on the set $\pi_0(R(\Omega_{MS}))$ we get that  for every $\Lambda'$ in $\Omega_{\Lambda}^0$ and $t$ in $\R^d$,
$$
\pi_1(\Lambda'-t) =  \pi_1(\Lambda')-t. 
$$

\subsubsection{Proof of \eqref{e:mtl}}

 Fix $\Lambda_1$  in  $R_{\pi_0}(\Omega_{MS})$. We prove that \eqref{e:mtl} holds. 
First, we assume that  $\Lambda_1$  is in $\pi_1(\Omega^0_{\Lambda}\cap \Xi^\Gamma)$. By Theorem \ref{t:parametmodel} if $\pi_0(\Lambda_1)=[(t,w)]$ then 
\begin{equation}\label{e:non-singular}
\curlywedge(w+\text{int}(W')) = \Lambda_1+t = \curlywedge(w+W').
\end{equation}
Observe that by definition of $\widetilde{\pi}$ for every  $\Lambda'$ in $\Omega^0_{\Lambda}\cap \Xi^\Gamma$ we have $\widetilde{\pi}(\Lambda')=[(0,\omega(\Lambda')]$. 
In addition, if $\Lambda'$ satisfies that $\pi_1(\Lambda')=\Lambda_1$ then using \eqref{eq: rel pi_0 pi_1}, we get
$$\pi_0(\Lambda_1)=\pi_0\circ\pi_1(\Lambda')=\widetilde{\pi}(\Lambda')=[(0,\omega(\Lambda'))].$$
Together with \eqref{e:non-singular} implies that  for every $\Lambda'\in\Omega^0_{\Lambda}\cap \Xi^\Gamma$ such that $\pi_1(\Lambda')=\Lambda_1$, we have 
\begin{equation}\label{e:Lambbda1}
\Lambda_1 = \{\gamma \in \Gamma \mid  s_{H'}(\gamma) \in \omega(\Lambda') + W'\}.
\end{equation}
By Proposition \ref{p:window}  and part (2) of Lemma \ref{l:omega}, we have
\begin{equation}\label{e:omega and the window}
-\omega(\Xi_\Lambda)=\theta([\Xi_\Lambda]_{\text{srp}})=W'.
\end{equation}
Since $\Lambda'\in\Xi^\Gamma$, and for every $\gamma\in\Lambda'$ we have  $\Lambda'-\gamma\in\Xi_{\Lambda}$, using part (1) of
 Lemma \ref{l:omega},  we get $\omega(\Lambda'-\gamma)=\omega(\Lambda')-s_{H'}(\gamma)$. Together with 
 \eqref{e:Lambbda1} and \eqref{e:omega and the window} this implies that for every $\gamma$ in $\Lambda'$ we have
\begin{multline*}
\omega(\Lambda'-\gamma) \in \omega(\Xi_{\Lambda}) \Longleftrightarrow s_{H'}(\gamma) \in \omega(\Lambda') -\omega(\Xi_{\Lambda})\\
\Longleftrightarrow s_{H'}(\gamma) \in \omega(\Lambda') +W'\Longleftrightarrow \gamma\in \Lambda_1.
\end{multline*}
Therefore, for every $\Lambda'\in\Omega^0_{\Lambda}\cap \Xi^\Gamma$ such that $\pi_1(\Lambda')=\Lambda_1$ we have 
\begin{equation}\label{e:in Lambda1}
\Lambda'\subseteq \Lambda_1.
\end{equation}

On the other hand, fix $\gamma$ in $\Lambda_1$. By   \eqref{e:Lambbda1} for every $\Lambda'\in\Omega^0_{\Lambda}\cap \Xi^\Gamma$ such that $\pi_1(\Lambda')=\Lambda_1$ we have  
$$s_{H'}(\gamma)\in \omega(\Lambda') + W' \Leftrightarrow  \omega(\Lambda')  \in  \omega(\Xi_{\Lambda}+\gamma).$$  
Thus, there is $\Lambda''$ in $\Xi_{\Lambda}+\gamma\subseteq \Xi^\Gamma$ such that $\omega(\Lambda'')=\omega(\Lambda')$. 
Then, $\Lambda''-\gamma$  is in $\Xi_{\Lambda}$ and thus, 
$\gamma$ is in $\Lambda''$. 
Therefore,
%$\pi_1(\Lambda'')=\pi_1(\Lambda')=\Lambda_1$ (FALTA VER QUE $\Lambda''$ ESTA EN $\Omega^0_\Lambda$) and
\begin{equation}\label{e:Lambda1 in 0}
\Lambda_1 \subseteq \bigcup_{\Lambda''\in \Omega^0_\Lambda \cap \Xi^\Gamma\text{ s.t. } \omega(\Lambda'') = \omega(\Lambda')} \Lambda''.
\end{equation}

Observe that  that for every  $\Lambda'\in\Omega^0_{\Lambda}\cap \Xi^\Gamma$ and every  $\Lambda'' \in \Xi^\Gamma$ such that 
$\omega(\Lambda')=\omega(\Lambda'')$ we have that 
$\tpi(\Lambda')=\tpi(\Lambda'')$ and thus,   $\Lambda''\in\Omega^0_{\Lambda}\cap \Xi^\Gamma$. In particular, 
$\pi_1(\Lambda')=\pi_1(\Lambda'')$, which together with \eqref{e:Lambda1 in 0}
implies 
\begin{equation}\label{e:Lambda1 in}
\Lambda_1 \subseteq \bigcup_{\pi_1(\Lambda'') = \Lambda_1} \Lambda''.
\end{equation}

Now, we prove that for every  $\Lambda'\in\Omega^0_{\Lambda}\cap \Xi^\Gamma$ and every  $\Lambda'' \in \Omega^0_\Lambda$ such that 
$\pi_1(\Lambda')=\pi_1(\Lambda'')$, we have that 
\begin{equation}\label{e:pi equivalences}
\Lambda''\in \Xi^\Gamma.
\end{equation}
First, observe that the definition of $\pi_1$ for all $\Lambda'$ and $\Lambda''$ in $\Omega^0_\Lambda$ we have that  
$\pi_1(\Lambda')=\pi_1(\Lambda'')  \Leftrightarrow  \widetilde{\pi}(\Lambda'') = \widetilde{\pi}(\Lambda')$. Now, let  $\Lambda'\in\Omega^0_{\Lambda}\cap \Xi^\Gamma$ and $\Lambda'' \in \Omega^0_\Lambda$  be such 
 that  $\widetilde{\pi}(\Lambda'') = \widetilde{\pi}(\Lambda')$. By definition of $\widetilde{\pi}$ this holds 
if and only if there exists $t$ in $\R^d$ such that $\Lambda''-t\in \Xi^\Gamma$ and $[(-t,\omega(\Lambda''-t))]=[(0,\omega(\Lambda'))]$, which is equivalent to 
the existence of $\gamma$ in $\Gamma$ such that 
$$(-t,\omega(\Lambda''-t)) - (0,\omega(\Lambda')) = (\gamma,s_{H'}(\gamma)).$$
Then, $-t = \gamma \in \Gamma$ and we get
$\Lambda''\subseteq \Gamma-\gamma=\Gamma$, which proves  \eqref{e:pi equivalences}. By \eqref{e:in Lambda1}, \eqref{e:Lambda1 in} and \eqref{e:pi equivalences} we conclude that 
\begin{equation*}
\Lambda_1 = \bigcup_{\pi_1(\Lambda'')=\Lambda_1} \Lambda'',
\end{equation*}
which is equivalent to
\begin{equation}\label{e:inclusions}
\Lambda_1 = \bigcup_{\Lambda'\in\; \widetilde{\pi}^{-1}(\pi_0(\Lambda_1))} \Lambda'.
\end{equation}

If $\Lambda_1$ is not in $\pi_1(\Omega_\Lambda^0\cap \Xi^\Gamma)$ then there is $t$ in $\R^d$ such that $\Lambda_1-t$ is in $\pi_1(\Omega_\Lambda^0\cap \Xi^\Gamma)$.
By \eqref{e:inclusions} we have that 
$$
\Lambda_1-t = \bigcup_{\tLambda\in\; \widetilde{\pi}^{-1}(\pi_0(\Lambda_1-t))} \tLambda.$$
Since $\tpi(\tLambda)= \pi_0(\Lambda_1-t)$ if and only if $\tpi(\tLambda-(-t))= \pi_0(\Lambda_1)$, we conclude that 
\begin{equation}\label{e:equality}
\Lambda_1 = \bigcup_{\tLambda\in\; \widetilde{\pi}^{-1}(\pi_0(\Lambda_1-t))} \tLambda - (-t) = \bigcup_{\Lambda'\in\; \widetilde{\pi}^{-1}(\pi_0(\Lambda_1))} \Lambda',
\end{equation}
which finishes the proof of \eqref{e:mtl}.

\subsubsection{$(\Omega_\Lambda,\R^d)$  almost automorphic extension of $(\Omega_{MS,\text{me}},\R^d)$}
Finally,  if $\tpi$ is an almost automorphic extension of $(\Omega_{MS,\text{me}},\R^d)$. 
By \cite[Lemma 4.1]{Ve}, we have that the set 
$$\tpi(R_{\tpi}(\Omega_\Lambda)) = \{x\in\Omega_{MS,\text{me}}\mid\widetilde{\pi}^{-1}(x)\;\textrm{is a singleton}\}$$
is a  residual set in $\Omega_{MS,\text{me}}$ and by \cite{Au16a}, the set $\pi_0(R_{\pi_0}(\Omega_{MS}))$ is also a residual set in 
$\Omega_{MS,\text{me}}$. Then, 
$$
\pi_0(R_{\pi_0}(\Omega_{MS}))\cap\tpi(R_{\tpi}(\Omega_\Lambda))
$$
is also a residual set in $\Omega_{MS,\text{me}}$. 
By \eqref{e:equality} for every $\Lambda_1$ in $R_{\pi_0}(\Omega_{MS})$ such that 
$\pi_0(\Lambda_1)\in \tpi(R_{\tpi}(\Omega_\Lambda))$ we have that $\Lambda_1$ is in $\Omega^0_\Lambda$,
which concludes the proof of the Main Technical Lemma'.

\bibliographystyle{alpha}
 %\bibliography{bibliomeyer.bib}
\newcommand{\etalchar}[1]{$^{#1}$}

\end{document}